\documentclass[11pt]{article}
\usepackage[dvips]{graphicx}
\usepackage{amsmath}
\usepackage{amsfonts}
\usepackage{amssymb}
\newtheorem{theorem}{Theorem}[subsection]

\newtheorem{conjecture}[theorem]{Conjecture}
\newtheorem{corollary}[theorem]{Corollary}

\newtheorem{example}[theorem]{Example}

\newtheorem{lemma}[theorem]{Lemma}

\newtheorem{proposition}[theorem]{Proposition}

\newenvironment{proof}[1][Proof]{\textbf{#1.} }{\ \rule{0.5em}{0.5em}}

\begin{document}

\title{A duality between $q$-multiplicities in tensor products and $q$-multiplicities
of weights for the root systems $B,C$ or $D$}
\author{C\'{e}dric Lecouvey\\lecouvey@math.unicaen.fr}
\date{}
\maketitle
\begin{abstract}
Starting from Jacobi-Trudi's type determinental expressions for the Schur
functions corresponding to types $B,C$ and $D,$ we define a natural
$q$-analogue of the multiplicity $[V(\lambda):M(\mu)]$ when $M(\mu)$ is a
tensor product of row or column shaped modules defined by $\mu$. We prove that
these $q$-multiplicities are equal to certain Kostka-Foulkes polynomials
related to the root systems $C$ or $D$.\ Finally we derive formulas expressing
the associated multiplicities in terms of Kostka numbers.
\end{abstract}

\section{Introduction}

Given two partitions $\lambda$ and $\mu$ of length $m,$ the Kostka number
$K_{\lambda,\mu}^{A_{m}}$ is equal to the dimension of the weight space $\mu$
in the finite dimensional irreducible $sl_{m}$-module $V^{A_{m}}(\lambda)$ of
highest weight $\lambda.\;$The Schur duality is a classical result in
representation theory establishing that $K_{\lambda,\mu}^{A_{m}}$ is also
equal to the multiplicity of $V(\lambda)$ in the tensor products%
\[
V(\mu_{1}\Lambda_{1})\otimes\cdot\cdot\cdot\otimes V(\mu_{m}\Lambda_{1})\text{
and }V(\Lambda_{\mu_{1}^{\prime}})\otimes\cdot\cdot\cdot\otimes V(\Lambda
_{\mu_{n}^{\prime}})
\]
where $\mu^{\prime}=(\mu_{1}^{\prime},...,\mu_{n}^{\prime})$ is the conjugate
partition of $\mu$ and the $\Lambda_{i}$'s $i=1,...,m-1$ are the fundamental
weights of $sl_{m}.$ Another way to define $K_{\lambda,\mu}^{A_{m}}$ is to use
the Jacobi-Trudi identity which gives a determinantal expression of the Schur
function $s_{\mu}=\mathrm{char}(V(\mu))$ in terms of the characters
$h_{k}=\mathrm{char}(V(k\Lambda_{1}))$ of the $k$-th symmetric power
representation.\ This formula can be rewritten%
\begin{equation}
s_{\mu}=\prod_{1\leq i<j\leq m}(1-R_{i,j})h_{\mu}\label{s-en-h}%
\end{equation}
where $h_{\mu}=h_{\mu_{1}}\cdot\cdot\cdot\cdot h_{\mu_{m}}$ and the $R_{i,j}$
are the raising operators (see \ref{raising}).\ Then one can prove that it
makes sense to write%
\begin{equation}
h_{\mu}=\prod_{1\leq i<j\leq m}(1-R_{i,j})^{-1}s_{\mu}\label{h_en_s}%
\end{equation}
which gives the decomposition of $h_{\mu}$ on the basis of Schur
functions.\ From this decomposition we derive the following expression for
$K_{\lambda,\mu}^{A_{m}}$:%
\begin{equation}
K_{\lambda,\mu}^{A_{m}}=\sum_{\sigma\in\mathcal{S}_{m}}(-1)^{l(\sigma
)}\mathcal{P}^{A_{m}}(\sigma(\lambda+\rho)-(\mu+\rho))\label{koskaA}%
\end{equation}
where $\mathcal{S}_{m}$ is the symmetric group of order $m$ and $\mathcal{P}%
^{A_{m}}$ the ordinary Kostant's partition function defined from the equality:%
\[
\prod_{\alpha\text{ positive root}}\dfrac{1}{(1-x^{\alpha})}=\sum_{\beta
}\mathcal{P}^{A_{m}}(\beta)x^{\beta}%
\]
with $\beta$ running on the set of nonnegative integral combinations of
positive roots of $sl_{m}$.

\noindent There exists a $q$-analogue $K_{\lambda,\mu}^{A_{m}}(q)$ of
$K_{\lambda,\mu}^{A_{m}}$ obtained by replacing the ordinary Kostant's
partition function $\mathcal{P}^{A_{m}}$ by its $q$-analogue $\mathcal{P}%
_{q}^{A_{m}}$ satisfying%
\[
\prod_{\alpha\text{ positive root}}\dfrac{1}{(1-qx^{\alpha})}=\sum_{\beta
}\mathcal{P}_{q}^{A_{m}}(\beta)x^{\beta}.
\]
So we have%
\begin{equation}
K_{\lambda,\mu}^{A_{m}}(q)=\sum_{\sigma\in\mathcal{S}_{m}}(-1)^{l(\sigma
)}\mathcal{P}_{q}^{A_{m}}(\sigma(\lambda+\rho)-(\mu+\rho)) \label{q-kostkaA}%
\end{equation}
which is a polynomial in $q$ with nonnegative integer coefficients
\cite{LSc1}, \cite{Lu}. In \cite{NY}, Nakayashiki and Yamada have shown that
$K_{\lambda,\mu}^{A_{m}}(q)$ can also be computed from the combinatorial $R$
matrix corresponding to Kashiwara's crystals associated to some $U_{q}%
(\widehat{sl_{m}})$-modules.

For $g=so_{2m+1},sp_{2m}$ or $so_{2m}$ there also exist expressions similar to
(\ref{koskaA}) for the multiplicities $K_{\lambda,\mu}^{g}$ of the weight
$\mu$ in the finite dimensional irreducible module $V(\lambda)$ but a so
simple duality as for $sl_{m}$ does not exist although it is possible to
obtain certain duality results between multiplicities of weights and tensor
product multiplicities of representations by using duals pairs of algebraic
groups (see \cite{HO}). This implies that the quantifications of weight
multiplicities and tensor product multiplicities can not coincide for the root
systems $B_{m},C_{m}$ and $D_{m}.\;$The Kostka-Foulkes polynomials
$K_{\lambda,\mu}^{g}(q)$ are the $q$-analogues of $K_{\lambda,\mu}^{g}$
defined as in (\ref{q-kostkaA}) by quantifying the partition function
corresponding to the root system associated to $g$ (see \ref{sub_sec_kost}).
In \cite{Ok}, Hatayama, Kuniba, Okado and Takagi have introduced for type
$C_{m}$ a quantification $X_{\lambda,\mu}^{C_{m}}(q)$ of the multiplicity of
$V(\lambda)$ in the tensor product%
\[
W(\mu_{1}\Lambda_{1})\otimes\cdot\cdot\cdot\otimes W(\mu_{m}\Lambda_{1})
\]
where for any $i=1,...,m,$
\[
W(\mu_{i}\Lambda_{1})=V(\mu_{i}\Lambda_{1})\oplus V((\mu_{i}-2)\Lambda
_{1})\oplus\cdot\cdot\cdot\oplus V((\mu_{i}\mathrm{mod}2)\Lambda_{1}).
\]
This quantification is based on the determination of the combinatorial $R$
matrix of some $U_{q}^{\prime}(g)$-crystals in the spirit of \cite{NY}. Note
that there also exit $q$-multiplicities for the $sp_{2}$-module $V(\lambda)$
in a tensor product
\[
V(\Lambda_{1})^{\otimes k}\otimes V(\Lambda_{2})^{\otimes l}%
\]
where $k,l$ are positive integers obtained by Yamada \cite{Y}.

In this paper we first use Jacobi-Trudi's type determinantal expressions for
the Schur functions associated to $g$ to introduce $q$-analogues of the
multiplicity of $V(\lambda)$ in the tensor products%
\begin{align*}
\mathrm{(i)} &  :\frak{h}(\mu)=V(\mu_{1}\Lambda_{1})\otimes\cdot\cdot
\cdot\otimes V(\mu_{m}\Lambda_{1}),\frak{H}(\mu)=W(\mu_{1}\Lambda_{1}%
)\otimes\cdot\cdot\cdot\otimes W(\mu_{m}\Lambda_{1})\\
\mathrm{(ii)} &  :\frak{e}(\mu)=V(\Lambda_{\mu_{1}^{\prime}})\otimes\cdot
\cdot\cdot\otimes V(\Lambda_{\mu_{n}^{\prime}}),\frak{E}(\mu)=W(\Lambda
_{\mu_{1}^{\prime}})\otimes\cdot\cdot\cdot\otimes W(\Lambda_{\mu_{n}^{\prime}%
})\text{ with }m\geq\left|  \mu\right|
\end{align*}
where
\[
\left\{
\begin{array}
[c]{c}%
W(\mu_{i}\Lambda_{1})=V(\mu_{i}\Lambda_{1})\oplus V((\mu_{i}-2)\Lambda
_{1})\oplus\cdot\cdot\cdot\oplus V((\mu_{i}\mathrm{mod}2)\Lambda_{1})\\
W(\Lambda_{k})=V(\Lambda_{k})\oplus V(\Lambda_{k-2})\oplus\cdot\cdot
\cdot\oplus V(\Lambda_{k\operatorname{mod}2})
\end{array}
\right.  \text{ }.
\]
With the condition $m\geq\left|  \mu\right|  $ for $\mathrm{(ii),}$ these
multiplicities are independent of the root system considered.\ When $q=1,$ we
recover a remarkable property already used by Koike and Terada in \cite{K}.
Next we prove that these $q$-multiplicities are in fact equal to
Kostka-Foulkes polynomials associated to the root systems of types $C$ and
$D.$ It is possible to extend the definition (\ref{q-kostkaA}) of the
Kostka-Foulkes polynomials associated to the root system $A_{m}$ by replacing
$\mu$ by $\gamma\in\mathbb{N}^{m}$ where $\gamma$ is not a partition.\ In this
case $K_{\lambda,\gamma}^{A_{m}}(q)$ may have nonnegative coefficients but
$K_{\lambda,\gamma}^{A_{m}}(1)$ is equal to the dimension of the weight space
$\gamma$ in $V^{A_{m}}(\lambda).$ We obtained simple expressions of the
$q$-multiplicities defined above in terms of the polynomials $K_{\lambda
,\gamma}^{A_{m}}(q).\;$By specializing at $q=1$ we derive formulas to compute
the related multiplicities from the Kostka numbers.

In section $1$ we recall background on the root systems $B_{m},C_{m}$ and
$D_{m}$ and the corresponding Kostka-Foulkes polynomials.\ We review in
section $2$ the determinantal identities for Schur functions that we need in
the sequel and we introduce the formalism suggested in \cite{LLTD} to prove
the expressions of Schur functions in terms of raising and lowering operators
implicitly contain in \cite{ram}.\ Thank to this formalism we are able to
obtain expressions for multiplicities similar to (\ref{koskaA}).\ We quantify
these multiplicities to obtain the desired $q$-analogues in section $3.\;$We
prove in Section $4$ two duality theorems between our $q$-analogues and
certain Kostka-Foulkes polynomials of types $C$ and $D.$ Finally we establish
formulas expressing the associated multiplicities in terms of Kostka numbers.

\bigskip

\noindent\textbf{Notation: }In the sequel we frequently define similar objects
for the root systems $B_{n}$ $C_{n}$ and $D_{n}$. When they are related to
type $B_{n}$ (resp. $C_{n},D_{n}$), we implicitly attach to them the label $B$
(resp. the labels $C,D$). To avoid cumbersome repetitions, we sometimes omit
the labels $B,C$ and $D$ when our definitions or statements are identical for
the three root systems.

\bigskip

\noindent\textbf{Note: }\textit{While writing this work, I have been informed
that Shimozono and Zabrocki \cite{SZ} have introduced independently and by
using creating operators essentially the same tensor power
multiplicities.\ Thanks to this formalism they recover in particular
Jacobi-Trudi's type determinantal expressions of the Schur functions
associated to the root systems }$B,C$ \textit{and }$D$ \textit{which
constitute the starting point of this article.}

\section{Background on the root systems $B_{m},C_{m}$ and $D_{m}$}

\subsection{Convention for the positive roots}

Consider an integer $m\geq1.$ The weight lattice for the root system $C_{m}$
(resp. $B_{m}$ and $D_{m})$ can be identified with $P_{C_{m}}=\mathbb{Z}^{m}$
(resp. $P_{B_{m}}=P_{D_{m}}\left(  \dfrac{\mathbb{Z}}{2}\right)  ^{m})$
equipped with the orthonormal basis $\varepsilon_{i},$ $i=1,...,m$.\ We take
for the simple roots%
\begin{equation}
\left\{
\begin{tabular}
[c]{l}%
$\alpha_{m}^{B_{m}}=\varepsilon_{m}\text{ and }\alpha_{i}^{B_{m}}%
=\varepsilon_{i}-\varepsilon_{i+1}\text{, }i=1,...,m-1\text{ for the root
system }B_{m}$\\
$\alpha_{m}^{C_{m}}=2\varepsilon_{m}\text{ and }\alpha_{i}^{C_{m}}%
=\varepsilon_{i}-\varepsilon_{i+1}\text{, }i=1,...,m-1\text{ for the root
system }C_{m}$\\
$\alpha_{m}^{D_{m}}=\varepsilon_{m}+\varepsilon_{m-1}\text{ and }\alpha
_{i}^{D_{m}}=\varepsilon_{i}-\varepsilon_{i+1}\text{, }i=1,...,m-1\text{ for
the root system }D_{m}$%
\end{tabular}
\right.  . \label{simple_roots}%
\end{equation}
Then the set of positive roots are%
\[
\left\{
\begin{tabular}
[c]{l}%
$R_{B_{m}}^{+}=\{\varepsilon_{i}-\varepsilon_{j},\varepsilon_{i}%
+\varepsilon_{j}\text{ with }1\leq i<j\leq m\}\cup\{\varepsilon_{i}\text{ with
}1\leq i\leq m\}\text{ for the root system }B_{m}$\\
$R_{C_{m}}^{+}=\{\varepsilon_{i}-\varepsilon_{j},\varepsilon_{i}%
+\varepsilon_{j}\text{ with }1\leq i<j\leq m\}\cup\{2\varepsilon_{i}\text{
with }1\leq i\leq m\}\text{ for the root system }C_{m}$\\
$R_{D_{m}}^{+}=\{\varepsilon_{i}-\varepsilon_{j},\varepsilon_{i}%
+\varepsilon_{j}\text{ with }1\leq i<j\leq m\}\text{ for the root system
}D_{m}$%
\end{tabular}
\right.  .
\]
Denote respectively by $P_{B_{m}}^{+},P_{C_{m}}^{+}$ and $P_{D_{m}}^{+}$the
sets of dominant weights of $so_{2m+1},sp_{2m}$ and $so_{2m}.$

\noindent Let $\lambda=(\lambda_{1},...,\lambda_{m})$ be a partition with $m$
parts. We will classically identify $\lambda$ with the dominant weight
$\sum_{i=1}^{m}\lambda_{i}\varepsilon_{i}.$ Note that there exists dominant
weights associated to the orthogonal root systems whose coordinates on the
basis $\varepsilon_{i},$ $i=1,...,m$ are not positive integers (hence which
can not be related to a partition). For each root system of type $B_{m},C_{m}$
or $D_{m},$ the set of weights having nonnegative integer coordinates on the
basis $\varepsilon_{1},...,\varepsilon_{m}$ can be identify with the set
$\pi_{m}^{+}$ of partitions of length $m.$\ For any partition $\lambda,$ the
weights of the finite dimensional $so_{2m+1},sp_{2m}$ or $so_{2m}$-module of
highest weight $\lambda$ are all in $\pi_{m}=\mathbb{Z}^{m}.\;$For any
$\alpha\in\pi_{m}$ we write $\left|  \alpha\right|  =\alpha_{1}+\cdot
\cdot\cdot\alpha_{m}$.

\noindent The conjugate partition of the partition $\lambda$ is denoted
$\lambda^{\prime}$ as usual.\ Consider $\lambda,\mu$ two partitions of length
$m$ and set $n=\mathrm{max}(\lambda_{1},\mu_{1})$.\ Then by adding to
$\lambda^{\prime}$ and $\mu^{\prime}$ the required numbers of parts $0$ we
will consider them as partitions of length $n.$

\bigskip

\noindent The Weyl group $W_{B_{m}}=W_{C_{m}}$ of $so_{2m+1}$ and $sp_{2m}$
can be regarded as the sub-group of the permutation group of $\{\overline
{m},...,\overline{2},\overline{1},1,2,...,m\}$\ generated by $s_{i}%
=(i,i+1)(\overline{i},\overline{i+1}),$ $i=1,...,m-1$ and $s_{m}%
=(m,\overline{m})$ where for $a\neq b$ $(a,b)$ is the simple transposition
which switches $a$ and $b.$ We denote by $l_{B}$ the length function
corresponding to the set of generators $s_{i},$ $i=1,...m.$

\noindent The Weyl group $W_{D_{m}}$ of $so_{2m}$ can be regarded as the sub
group of the permutation group of $\{\overline{m},...,\overline{2}%
,\overline{1},1,2,...,m\}$\ generated by $s_{i}=(i,i+1)(\overline{i}%
,\overline{i+1}),$ $i=1,...,m-1$ and $s_{m}^{\prime}=(m,\overline
{m-1})(m-1,\overline{m})$. We denote by $l_{D}$ the length function
corresponding to the set of generators $s_{m}^{\prime}$ and $s_{i},$ $i=1,...m-1.$

\noindent Note that $W_{D_{m}}\subset W_{B_{m}}$ and any $w\in W_{B_{m}}$
verifies $w(\overline{i})=\overline{w(i)}$ for $i\in\{1,...,m\}.$ The action
of $w$ on $\beta=(\beta_{1},...,\beta_{m})\in P_{m}$ is given by%
\[
w\cdot(\beta_{1},...,\beta_{m})=(\beta_{1}^{w},...,\beta_{m}^{w})
\]
where $\beta_{i}^{w}=\beta_{w(i)}$ if $\sigma(i)\in\{1,...,m\}$ and $\beta
_{i}^{w}=-\beta_{w(\overline{i})}$ otherwise.

\noindent The half sums $\rho_{B_{m}},\rho_{C_{m}}$ and $\rho_{D_{m}}$ of the
positive roots associated to each root system $B_{m},C_{m}$ and $D_{m}$
verify:%
\[
\rho_{B_{m}}=(m-\dfrac{1}{2},m-\dfrac{3}{2},...,\dfrac{1}{2}),\rho_{C_{m}%
}=(m,m-1,...,1)\text{ and }\rho_{B_{m}}=(m-1,m-2,...,0).
\]
In the sequel we identify the symmetric group $\mathcal{S}_{m}$ with the sub
group of $W_{B_{m}}$ or $W_{D_{m}}$ generated by the $s_{i}$'s, $i=1,...,m-1.$

\subsection{Schur functions and Kostka-Foulkes polynomials\label{sub_sec_kost}}

We now briefly review the notions of Schur functions and Kostka-Foulkes
polynomials associated to the roots systems $B_{m},C_{m}$ and $D_{m}$ and
refer the reader to \cite{NR} for more details.\ For any weight $\beta
=(\beta_{1},...,\beta_{m})\in\pi_{m}$ we set $x^{\beta}=x_{1}^{\beta_{1}}%
\cdot\cdot\cdot x_{m}^{\beta_{m}}$ where $x_{1},...,x_{m}$ are fixed
indeterminates. We set
\[
a_{\beta}=\sum_{w\in W_{B_{m}}}(-1)^{l(\sigma)}(w\cdot x^{\beta})
\]
where $w\cdot x^{\mu}=x^{w(\mu)}.$ The Schur function $s_{\beta}$ is defined
by
\[
s_{\beta}^{B_{m}}=\dfrac{a_{\beta+\rho_{B_{m}}}}{a_{\rho_{B_{m}}}}.
\]
When $\nu\in\pi_{m}^{+},$ $s_{\nu}^{B_{m}}$ is the Weyl character of $V(\nu)$
the finite dimensional irreducible module with highest weight $\nu.$ For any
$w\in W_{B_{m}},$ the dot action of $w$ on $\beta\in\pi_{m}$ is defined by
\[
w\circ\beta=w\cdot(\beta+\rho_{B_{m}})-\rho_{B_{m}}.
\]
We have the following straightening law for the Schur functions. For any
$\beta\in\pi_{m}$, $s_{\beta}^{B_{m}}=0$ or there exists a unique $\nu\in
\pi_{m}^{+}$ such that $s_{\beta}^{B_{m}}=(-1)^{l(w)}s_{\nu}^{B_{m}}$ with
$w\in W_{B_{m}}$ and $\nu=w\circ\beta.$ Set $\mathbb{K}=\mathbb{Z}[q,q^{-1}]$
and write $\mathbb{K}[\pi_{m}]$ for the $\mathbb{K}$-module generated by the
$x^{\beta}$, $\beta\in\pi_{m}.$ Set $\mathcal{C}_{B_{m}}=\mathbb{K}[\pi
_{m}]^{W_{B_{m}}}=\{f\in\mathbb{K}[\pi_{m}],$ $w\cdot f=f$ for any $w\in
W_{B_{m}}\}.$ Then $\{s_{\nu}^{B_{m}}\},\nu\in\pi_{m}^{+}$ is a basis of
$\mathbb{K}[\pi_{m}]^{W_{B_{m}}}.$

\noindent We define $s_{\beta}^{C_{m}}$ and $s_{\beta}^{D_{m}}$ belonging to
$\mathcal{C}_{C_{m}}=\mathcal{C}_{B_{m}}$ and $\mathcal{C}_{D_{m}}$ in the
same way when $\beta\in\pi_{m}^{+}$ or $\beta\in\pi_{m}^{+}$ and we obtain
similarly that $\{s_{\nu}^{C_{m}}\},\nu\in\pi_{m}^{+}$ and $\{s_{\nu}^{D_{m}%
}\},\nu\in\pi_{m}^{+}$ are respectively bases of $\mathcal{C}_{C_{m}}$ and
$\mathcal{C}_{D_{m}}.$

\bigskip

\noindent The $q$-analogue $\mathcal{P}_{q}^{B_{m}}$ of the Kostant's
partition function corresponding to the root system $B_{m}$ is defined by the
equality%
\[
\prod_{\alpha\in R_{B_{m}}^{+}}\dfrac{1}{1-qx^{\alpha}}=\sum_{\beta\in\pi_{m}%
}\mathcal{P}_{q}^{B_{m}}(\beta)x^{\beta}.
\]
Note that $\mathcal{P}_{q}^{B_{m}}(\beta)=0$ if $\beta$ is not a linear
combination of positive roots of $R_{B_{m}}^{+}$with nonnegative
coefficients.\ We define similarly $\mathcal{P}_{q}^{C_{m}}$and $\mathcal{P}%
_{q}^{D_{m}}$the $q$-Kostant's partition functions respectively associated to
the root systems $C_{m}$ and $D_{m}$. Given $\lambda$ and $\mu$ two partitions
of length $m,$ the Kostka-Foulkes polynomials of types $B_{m},C_{m}$ and
$D_{m}$ are then respectively defined by%
\begin{align*}
K_{\lambda,\mu}^{B_{m}}(q) &  =\sum_{\sigma\in W_{B_{m}}}(-1)^{l(\sigma
)}\mathcal{P}_{q}^{B_{m}}(\sigma(\lambda+\rho_{B_{m}})-(\mu+\rho_{B_{m}})),\\
K_{\lambda,\mu}^{C_{m}}(q) &  =\sum_{\sigma\in W_{C_{m}}}(-1)^{l(\sigma
)}\mathcal{P}_{q}^{C_{m}}(\sigma(\lambda+\rho_{C_{m}})-(\mu+\rho_{C_{m}})),\\
K_{\lambda,\mu}^{D_{m}}(q) &  =\sum_{\sigma\in W_{D_{m}}}(-1)^{l(\sigma
)}\mathcal{P}_{q}^{D_{m}}(\sigma(\lambda+\rho_{D_{m}})-(\mu+\rho_{D_{m}})).
\end{align*}

\noindent\textbf{Remarks: }

\noindent$\mathrm{(i):}$ We have $K_{\lambda,\mu}(q)=0$ when $\left|
\lambda\right|  <\left|  \mu\right|  .$

\noindent$\mathrm{(ii):}$ When $\left|  \lambda\right|  =\left|  \mu\right|
,$ $K_{\lambda,\mu}^{B_{m}}(q)=K_{\lambda,\mu}^{C_{m}}(q)=K_{\lambda,\mu
}^{D_{m}}(q)=K_{\lambda,\mu}^{A_{m}}(q)$ that is, the Kostka-Foulkes
polynomials associated to the root systems $B_{m},C_{m}$ and $D_{m}$ are
Kostka-Foulkes polynomials associated to the root system $A_{m}.$

\section{Determinantal identities and multiplicities of representations}

\subsection{Determinantal identities for Schur functions}

Consider $k\in\mathbb{Z}$. When $k$ is a nonnegative integer, write
$(k)_{m}=(k,0,...,0)$ for the partition of length $m$ with a unique non-zero
part equal to $k$.\ Then set
\[
h_{k}^{B_{m}}=s_{(k)_{m}}^{B_{m}},h_{k}^{C_{m}}=s_{(k)_{m}}^{C_{m}}%
,h_{k}^{D_{m}}=s_{(k)_{m}}^{D_{m}}%
\]
and%
\begin{gather*}
H_{k}^{B_{m}}=h_{k}^{B_{m}}+h_{k-2}^{B_{m}}+\cdot\cdot\cdot
+h_{k\operatorname{mod}2}^{B_{m}},H_{k}^{C_{m}}=h_{k}^{C_{m}}+h_{k-2}^{C_{m}%
}+\cdot\cdot\cdot+h_{k\operatorname{mod}2}^{B_{m}},\\
H_{k}^{D_{m}}=h_{k}^{D_{m}}+h_{k-2}^{D_{m}}+\cdot\cdot\cdot
+h_{k\operatorname{mod}2}^{D_{m}}.
\end{gather*}
When $k$ is a negative integer we set $h_{k}^{B_{m}}=h_{k}^{C_{m}}%
=h_{k}^{D_{m}}=0$ and $H_{k}^{B_{m}}=H_{k}^{C_{m}}=H_{k}^{D_{m}}=0.$ For any
$\alpha=(\alpha_{1},...,\alpha_{m})\in\mathbb{Z}^{m}$ define
\[
u_{\alpha}^{B_{m}}=\det\left(
\begin{array}
[c]{cccc}%
h_{\alpha_{1}}^{B_{m}} & h_{\alpha_{1}+1}^{B_{m}}+h_{\alpha_{1}-1}^{B_{m}} &
\cdot\cdot\cdot\cdot\cdot\cdot\cdot\cdot\cdot\cdot &  h_{\alpha_{1}%
+m-1}^{B_{m}}+h_{\alpha_{1}-m+1}^{B_{m}}\\
h_{\alpha_{2}-1}^{B_{m}} & h_{\alpha_{2}}^{B_{m}}+h_{\alpha_{2}-2}^{B_{m}} &
\cdot\cdot\cdot\cdot\cdot\cdot\cdot\cdot\cdot\cdot &  h_{\alpha_{2}%
+m-2}^{B_{m}}+h_{\alpha_{2}-m}^{B_{m}}\\
\cdot & \cdot & \cdot\cdot\cdot\cdot\cdot\cdot\cdot\cdot\cdot\cdot & \cdot\\
\cdot & \cdot & \cdot\cdot\cdot\cdot\cdot\cdot\cdot\cdot\cdot\cdot & \cdot\\
h_{\alpha_{m}-m+1}^{B_{m}} & h_{\alpha_{m}-m+2}^{B_{m}}+h_{\alpha_{m}%
-m}^{B_{m}} & \cdot\cdot\cdot\cdot\cdot\cdot\cdot\cdot\cdot\cdot &
h_{\alpha_{m}}^{B_{m}}+h_{\alpha_{m}-2m+2}^{B_{m}}%
\end{array}
\right)  .
\]
By using the equalities $h_{k}^{B_{m}}=H_{k}^{B_{m}}-H_{k-2}^{B_{m}}$ and
simple computations on determinants we have also%

\begin{equation}
u_{\alpha}^{B_{m}}=\det\left(
\begin{array}
[c]{cccc}%
H_{\alpha_{_{1}}}^{B_{m}}-H_{\alpha_{1}-2}^{B_{m}} & H_{\alpha_{1}+1}^{B_{m}%
}-H_{\alpha_{1}-1}^{B_{m}} & \cdot\cdot\cdot\cdot\cdot\cdot\cdot\cdot
\cdot\cdot &  H_{\alpha_{1}+m-1}^{B_{m}}-H_{\alpha_{1}-m-1}^{B_{m}}\\
H_{\alpha_{2}-1}^{B_{m}}-H_{\alpha_{2}-3}^{B_{m}} & H_{\alpha_{2}}^{B_{m}%
}-H_{\alpha_{2}-4}^{B_{m}} & \cdot\cdot\cdot\cdot\cdot\cdot\cdot\cdot
\cdot\cdot &  H_{\alpha_{2}+m-2}^{B_{m}}-H_{\alpha_{2}-m-2}^{B_{m}}\\
\cdot & \cdot & \cdot\cdot\cdot\cdot\cdot\cdot\cdot\cdot\cdot\cdot & \cdot\\
\cdot & \cdot & \cdot\cdot\cdot\cdot\cdot\cdot\cdot\cdot\cdot\cdot & \cdot\\
H_{\alpha_{m}-m+1}^{B_{m}}-H_{\alpha_{m}-m-1}^{B_{m}} & H_{\alpha_{m}%
-m+2}^{B_{m}}-H_{\alpha_{m}-m-2}^{B_{m}} & \cdot\cdot\cdot\cdot\cdot\cdot
\cdot\cdot\cdot\cdot &  H_{\alpha_{m}}^{B_{m}}+H_{\alpha_{m}-2m-2}^{B_{m}}%
\end{array}
\right)  . \label{u_H}%
\end{equation}
We define $u_{\alpha}^{C_{m}}$ and $u_{\alpha}^{D_{m}}$similarly by replacing
$h_{k}^{B_{m}}$ respectively by $h_{k}^{C_{m}}$ and $h_{k}^{D_{m}}$.

Consider $p$ and $m$ two integers such that $m\geq1.$ When $p$ is nonnegative
and $m\geq p$, write $(1^{p})_{m}=(1,...,1,0,...,0)$ for the partition of
length $m$ having $p$ non zero parts equal to $1.$ We set
\[
\left\{
\begin{array}
[c]{l}%
e_{p}^{B_{m}}=s_{(1^{p})_{m}}^{B_{m}},e_{p}^{C_{m}}=s_{(1^{p})_{m}}^{C_{m}%
},e_{p}^{D_{m}}=s_{(1^{p})_{m}}^{D_{m}}\text{ if }0\leq p\leq m\\
e_{p}^{B_{m}}=e_{2p-m}^{B_{m}},e_{p}^{C_{m}}=e_{2p-m}^{C_{m}},e_{k}^{D_{m}%
}=e_{2p-m}^{D_{m}}\text{ if }m+1\leq p\leq2m\\
e_{p}^{B_{m}}=e_{p}^{C_{m}}=e_{p}^{D_{m}}=0\text{ otherwise}%
\end{array}
\right.  .
\]
and%
\begin{gather*}
E_{k}^{B_{m}}=e_{k}^{B_{m}}+e_{k-2}^{B_{m}}+\cdot\cdot\cdot
+e_{k\operatorname{mod}2}^{B_{m}},H_{k}^{C_{m}}=e_{k}^{C_{m}}+e_{k-2}^{C_{m}%
}+\cdot\cdot\cdot+e_{k\operatorname{mod}2}^{B_{m}},\\
H_{k}^{D_{m}}=e_{k}^{D_{m}}+e_{k-2}^{D_{m}}+\cdot\cdot\cdot
+e_{k\operatorname{mod}2}^{D_{m}}.
\end{gather*}
For any $\beta=(\beta_{1},...,\beta_{m})\in\mathbb{Z}^{m}$ we define
\[
v_{\beta}^{B_{m}}=\det\left(
\begin{array}
[c]{cccc}%
e_{\beta_{1}}^{B_{m}} & e_{\beta_{1}+1}^{B_{m}}+e_{\beta_{1}-1}^{B_{m}} &
\cdot\cdot\cdot\cdot\cdot\cdot\cdot\cdot\cdot\cdot &  e_{\beta_{1}+m-1}%
^{B_{m}}+e_{\beta_{1}-m+1}^{B_{m}}\\
e_{\beta_{2}-1}^{B_{m}} & e_{\beta_{2}}^{B_{m}}+e_{\beta_{2}-2}^{B_{m}} &
\cdot\cdot\cdot\cdot\cdot\cdot\cdot\cdot\cdot\cdot &  e_{\beta_{2}+m-2}%
^{B_{m}}+e_{\beta_{2}-m}^{B_{m}}\\
\cdot & \cdot & \cdot\cdot\cdot\cdot\cdot\cdot\cdot\cdot\cdot\cdot & \cdot\\
\cdot & \cdot & \cdot\cdot\cdot\cdot\cdot\cdot\cdot\cdot\cdot\cdot & \cdot\\
e_{\beta_{m}-m+1}^{B_{m}} & e_{\beta_{m}-m+2}^{B_{m}}+e_{\beta_{m}-m}^{B_{m}}%
& \cdot\cdot\cdot\cdot\cdot\cdot\cdot\cdot\cdot\cdot &  e_{\beta_{m}}^{B_{m}%
}+e_{\beta_{m}-2m+2}^{B_{m}}%
\end{array}
\right)  .
\]
By using the equalities $e_{k}^{B_{m}}=E_{k}^{B_{m}}-E_{k-2}^{B_{m}}$ and
simple computations on determinants we have also%
\[
v_{\beta}^{B_{m}}=\det\left(
\begin{array}
[c]{cccc}%
E_{\beta_{_{1}}}^{B_{m}}-E_{\beta_{1}-2}^{B_{m}} & E_{\beta_{1}+1}^{B_{m}%
}-E_{\beta_{1}-1}^{B_{m}} & \cdot\cdot\cdot\cdot\cdot\cdot\cdot\cdot\cdot\cdot
&  E_{\beta_{1}+m-1}^{B_{m}}-E_{\beta_{1}-m-1}^{B_{m}}\\
E_{\beta_{2}-1}^{B_{m}}-E_{\beta_{2}-3}^{B_{m}} & E_{\beta_{2}}^{B_{m}%
}-E_{\beta_{2}-4}^{B_{m}} & \cdot\cdot\cdot\cdot\cdot\cdot\cdot\cdot\cdot\cdot
&  E_{\beta_{2}+m-2}^{B_{m}}-E_{\beta_{2}-m-2}^{B_{m}}\\
\cdot & \cdot & \cdot\cdot\cdot\cdot\cdot\cdot\cdot\cdot\cdot\cdot & \cdot\\
\cdot & \cdot & \cdot\cdot\cdot\cdot\cdot\cdot\cdot\cdot\cdot\cdot & \cdot\\
E_{\beta_{m}-m+1}^{B_{m}}-E_{\beta_{m}-m-1}^{B_{m}} & E_{\beta_{m}-m+2}%
^{B_{m}}-E_{\beta_{m}-m-2}^{B_{m}} & \cdot\cdot\cdot\cdot\cdot\cdot\cdot
\cdot\cdot\cdot &  E_{\beta_{m}}^{B_{m}}+E_{\beta_{m}-2m-2}^{B_{m}}%
\end{array}
\right)
\]
The determinants $v_{\beta}^{C_{m}},v_{\beta}^{D_{m}}$ are defined similarly.

\begin{proposition}
\label{prop_fult}(see\cite{FH}) Consider $\lambda$ a partition of length $m$
and suppose that $\lambda^{\prime}=(\lambda_{1}^{\prime},...,\lambda
_{n}^{\prime})$ is a partition of length $n$. Then $u_{\lambda}=s_{\lambda}$
and $v_{\lambda^{\prime}}=s_{\lambda}.$
\end{proposition}

\begin{lemma}
\label{straight}(straigthening law for $u_{\alpha}$ and $v_{\beta})$

\noindent Consider $\alpha\in\pi_{m}^{+}$ then
\[
u_{\alpha}=\left\{
\begin{tabular}
[c]{l}%
$(-1)^{l(\sigma)}u_{\lambda}$ if there exists $\sigma\in\mathcal{S}_{m}$ and
$\lambda\in\pi_{m}^{+}$ such that $\sigma\circ\alpha=\lambda$\\
$0$ otherwise
\end{tabular}
\right.  .
\]

\noindent Consider $\beta\in\pi_{n}^{+}$ then
\[
v_{\beta}=\left\{
\begin{tabular}
[c]{l}%
$(-1)^{l(\sigma)}v_{\nu}$ if there exists $\sigma\in\mathcal{S}_{n}$ and
$\nu\in\pi_{n}^{+}$ such that $\sigma\circ\alpha=\nu$\\
$0$ otherwise
\end{tabular}
\right.  .
\]
\end{lemma}

\begin{proof}
By commutting the rows $i$ and $i+1$ in the determinant (\ref{u_H}) we see
that $u_{s_{i}\circ\alpha}=-u_{\alpha}.\;$This implies that $u_{\sigma
\circ\alpha}=(-1)^{l(\sigma)}u_{\alpha}$ for any $\sigma\in\mathcal{S}_{m}%
.\;$Then it follows from the definition of the dot action that $u_{\alpha}=0$
or there exists $\gamma\in\pi_{m}$ and $\sigma\in\mathcal{S}_{m}$ such that
$\gamma_{1}\geq\cdot\cdot\cdot\geq\gamma_{m}$ and $\gamma=\sigma\circ
\alpha.\;$In this last case we have $u_{\alpha}=(-1)^{l(\sigma)}u_{\gamma}%
$.\ Now if there exists a negative $\gamma_{i}$, $u_{\gamma}=0$ since all the
$H_{k}$ which appear in the lowest row of (\ref{u_H}) are equal to $0.$ The
proof is similar for $v_{\beta}.$
\end{proof}

\subsection{Determinantal identities in terms of raising and lowering
operators\label{raising}}

Let $\mathcal{L}_{m}=\mathbb{K[}[x_{1},x_{1}^{-1},...,x_{m},x_{m}^{-1}]]$ be
the ring of formal series in the indeterminates $x_{1},x_{1}^{-1}%
,...,x_{m},x_{m}^{-1}.\;$We consider the two following determinants%
\begin{align*}
\delta_{m}(\alpha) &  =\det\left(
\begin{array}
[c]{cccc}%
x_{1}^{\alpha_{1}} & x_{1}^{\alpha_{1}+1}+x_{1}^{\alpha_{1}-1} & \cdot
\cdot\cdot\cdot\cdot\cdot\cdot\cdot\cdot\cdot &  x_{1}^{\alpha_{1}+m-1}%
+x_{1}^{\alpha_{1}-m+1}\\
x_{2}^{\alpha_{2}} & x_{2}^{\alpha_{2}}+x_{2}^{\alpha_{2}-2} & \cdot\cdot
\cdot\cdot\cdot\cdot\cdot\cdot\cdot\cdot &  x_{2}^{\alpha_{2}+m-2}%
+x_{2}^{\alpha_{2}-m}\\
\cdot & \cdot & \cdot\cdot\cdot\cdot\cdot\cdot\cdot\cdot\cdot\cdot & \cdot\\
\cdot & \cdot & \cdot\cdot\cdot\cdot\cdot\cdot\cdot\cdot\cdot\cdot & \cdot\\
x_{m}^{\alpha_{m}-m+1} & x_{m}^{\alpha_{m}-m+2}+x_{m}^{\alpha_{m}-m} &
\cdot\cdot\cdot\cdot\cdot\cdot\cdot\cdot\cdot\cdot &  x_{m}^{\alpha_{m}}%
+x_{m}^{\alpha_{m}-2m+2}%
\end{array}
\right)  \text{ and}\vspace{0.6cm}\\
\Delta_{m}(\alpha) &  =\det\left(
\begin{array}
[c]{cccc}%
x_{1}^{\alpha_{1}}-x_{1}^{\alpha_{1}-2} & x_{1}^{\alpha_{1}+1}-x_{1}%
^{\alpha_{1}-1} & \cdot\cdot\cdot\cdot\cdot\cdot\cdot\cdot\cdot\cdot &
x_{1}^{\alpha_{1}+m-1}-x_{1}^{\alpha_{1}-m-1}\\
x_{2}^{\alpha_{2}-1}-x_{2}^{\alpha_{2}-3} & x_{2}^{\alpha_{2}}-x_{2}%
^{\alpha_{2}-4} & \cdot\cdot\cdot\cdot\cdot\cdot\cdot\cdot\cdot\cdot &
x_{2}^{\alpha_{2}+m-2}-x_{2}^{\alpha_{2}-m-2}\\
\cdot & \cdot & \cdot\cdot\cdot\cdot\cdot\cdot\cdot\cdot\cdot\cdot & \cdot\\
\cdot & \cdot & \cdot\cdot\cdot\cdot\cdot\cdot\cdot\cdot\cdot\cdot & \cdot\\
x_{m}^{\alpha_{m}-m+1}+x_{m}^{\alpha_{m}-m-1} & x_{m}^{\alpha_{m-m+2}}%
-x_{m}^{\alpha_{m}-m} & \cdot\cdot\cdot\cdot\cdot\cdot\cdot\cdot\cdot\cdot &
x_{m}^{\alpha_{m}}-x_{m}^{\alpha_{m}-2m-2}%
\end{array}
\right)
\end{align*}
From a simple computation we derive the equalities:%
\begin{equation}
\delta_{m}(\alpha)=\prod_{1\leq i<j\leq m}(1-\frac{x_{i}}{x_{j}})\prod_{1\leq
r<s\leq m}(1-\frac{1}{x_{i}x_{j}})x^{\alpha}\text{ and }\Delta_{m}%
(\alpha)=\prod_{1\leq i<j\leq m}(1-\frac{x_{i}}{x_{j}})\prod_{1\leq r\leq
s\leq m}(1-\frac{1}{x_{i}x_{j}})x^{\alpha}.\label{delda_factor}%
\end{equation}
We set $h_{\alpha}=h_{\alpha_{1}}\cdot\cdot\cdot h_{\alpha_{m}},H_{\alpha
}=H_{\alpha_{1}}\cdot\cdot\cdot H_{\alpha_{m}},e_{\alpha}=e_{\alpha_{1}}%
\cdot\cdot\cdot e_{\alpha_{m}}$ and $E_{\alpha}=E_{\alpha_{1}}\cdot\cdot\cdot
E_{\alpha_{m}}.$

\noindent\textbf{Remarks}

\noindent$\mathrm{(i):}$ For any partition $\mu$ of length $m,$ $h_{\mu}$ is
the character of $\frak{h}(\mu)=V(\mu_{1}\Lambda_{1})\otimes\cdot\cdot
\cdot\otimes V(\mu_{m}\Lambda_{1})$ and $H_{\mu}$ is the character of
$\frak{H}(\mu)=W(\mu_{1}\Lambda_{1})\otimes\cdot\cdot\cdot\otimes W(\mu
_{m}\Lambda_{1})$ where for any $k\in\mathbb{N}$, $W(k_{1})=V(k\Lambda
_{1})\oplus V((k-2)\Lambda_{1})\oplus\cdot\cdot\cdot\oplus
V((k\operatorname{mod}2)\Lambda_{1}).$

\noindent$\mathrm{(ii):}$ For any partition $\mu$ of length $m$ such that
$\mu^{\prime}$ is of length $n,$ $e_{\mu^{\prime}}$ is the character of
$\frak{e}(\mu)=V(\Lambda_{\mu_{1}^{\prime}})\otimes\cdot\cdot\cdot\otimes
V(\Lambda_{\mu_{n}^{\prime}})$ and $E_{\mu^{\prime}}$ is the character of
$\frak{E}(\mu)=W(\Lambda_{\mu_{1}^{\prime}})\otimes\cdot\cdot\cdot\otimes
W(\Lambda_{\mu_{n}^{\prime}})$ where for any $k\in\mathbb{N}$ with $k\leq m$,
$W(\Lambda_{k})=V(\Lambda_{k})\oplus V(\Lambda_{k-2})\oplus\cdot\cdot
\cdot\oplus V(\Lambda_{k\operatorname{mod}2}).$

$\bigskip$

\noindent For the root system $B_{m}$ we introduce six linear maps
$\mathrm{h}_{B_{m}},\mathrm{H}_{B_{m}}$,$\mathrm{u}_{B_{m}}$ and
$\mathrm{e}_{B_{m}},\mathrm{E}_{B_{m}}$,$\mathrm{v}_{B_{m}}$ as follows:%
\begin{gather*}
\left\{
\begin{tabular}
[c]{c}%
$\mathrm{h}_{B_{m}}:\mathcal{L}_{m}\rightarrow\mathcal{C}_{B_{m}}$\\
$\ \ \ \ \ \ \ \ \ x^{\alpha}\mapsto h_{\alpha}^{B_{m}}$%
\end{tabular}
\right.  ,\left\{
\begin{tabular}
[c]{c}%
$\mathrm{H}_{B_{m}}:\mathcal{L}_{m}\rightarrow\mathcal{C}_{B_{m}}$\\
$x^{\alpha}\mapsto H_{\alpha}^{B_{m}}$%
\end{tabular}
\right.  ,\text{ }\left\{
\begin{tabular}
[c]{c}%
$\mathrm{u}_{B_{m}}:\mathcal{L}_{m}\rightarrow\mathcal{C}_{B_{m}}$\\
$x^{\alpha}\mapsto u_{\alpha}^{B_{m}}$%
\end{tabular}
\right.  \text{ and }\\
\left\{
\begin{tabular}
[c]{c}%
$\mathrm{e}_{B_{m}}:\mathcal{L}_{m}\rightarrow\mathcal{C}_{B_{m}}$\\
$x^{\alpha}\mapsto e_{\alpha}^{B_{m}}$%
\end{tabular}
\right.  ,\left\{
\begin{tabular}
[c]{c}%
$\mathrm{E}_{B_{m}}:\mathcal{L}_{m}\rightarrow\mathcal{C}_{B_{m}}$\\
$x^{\alpha}\mapsto E_{\alpha}^{B_{m}}$%
\end{tabular}
\right.  ,\text{ }\left\{
\begin{tabular}
[c]{c}%
$\mathrm{v}_{B_{m}}:\mathcal{L}_{m}\rightarrow\mathcal{C}_{B_{m}}$\\
$x^{\alpha}\mapsto v_{\alpha}^{B_{m}}$%
\end{tabular}
\right.  .
\end{gather*}
Note that these maps are not ring homomorphisms. For the roots systems $C_{m}$
and $D_{m}$ we define respectively the maps $\mathrm{h}_{C_{m}},\mathrm{H}%
_{C_{m}},\mathrm{u}_{C_{m}},\mathrm{e}_{C_{m}},\mathrm{E}_{C_{m}}%
,\mathrm{v}_{C_{m}}$ and $\mathrm{h}_{D_{m}},\mathrm{H}_{D_{m}},\mathrm{u}%
_{D_{m}},\mathrm{e}_{D_{m}},\mathrm{E}_{D_{m}},\mathrm{v}_{D_{m}}$ similarly.

\noindent Let $\mathrm{\omega}_{m}$ and $\mathrm{\Omega}_{m}$ be the
endomorphisms of $\mathcal{L}_{m}$ corresponding respectively to the
multiplication by
\[
\phi_{m}=\prod_{1\leq i<j\leq m}(1-\frac{x_{i}}{x_{j}})\prod_{1\leq r<s\leq
m}(1-\frac{1}{x_{i}x_{j}})\text{ and }\Phi_{m}=\prod_{1\leq i<j\leq m}%
(1-\frac{x_{i}}{x_{j}})\prod_{1\leq r\leq s\leq m}(1-\frac{1}{x_{i}x_{j}}).
\]
Since $\phi_{m}^{-1}$ and $\Phi_{m}^{-1}$ belong to $\mathcal{L}_{m},$
$\mathrm{\omega}_{m}$ and $\mathrm{\Omega}_{m}$ are the automorphisms of
$\mathcal{L}_{m}$ corresponding to the multiplication by $\phi_{m}^{-1}$ and
$\Phi_{m}^{-1}.$

\begin{proposition}
\label{prop_s-h}We have

\begin{enumerate}
\item $\mathrm{u}_{m}=\mathrm{h}_{m}\cdot\mathrm{\omega}_{m}$ and
$\mathrm{u}_{m}=\mathrm{H}_{m}\cdot\mathrm{\Omega}_{m},$

\item $\mathrm{v}_{m}=\mathrm{e}_{m}\cdot\mathrm{\omega}_{m}$ and
$\mathrm{v}_{m}=\mathrm{E}_{m}\cdot\mathrm{\Omega}_{m}.$
\end{enumerate}
\end{proposition}

\begin{proof}
$1:$ We have seen that $\mathrm{h}_{m}$ is not a
ring-homomorphism.\ Nevertheless we have by definition of the $h_{\alpha}$%
\[
\mathrm{h}_{m}(x^{\alpha})=\mathrm{h}_{m}(x_{1}^{\alpha_{1}})\cdot\cdot
\cdot\mathrm{h}_{m}(x_{m}^{\alpha_{m}})=h_{\alpha_{1}}\cdot\cdot\cdot
h_{\alpha_{m}}.
\]
More generally if $P_{1},...,P_{m}$ are polynomials respectively in the
indeterminates $x_{1},...,x_{m}$, we have%
\[
\mathrm{h}_{m}(P_{1}(x_{1})\cdot\cdot\cdot P_{m}(x_{m}))=\mathrm{h}_{m}%
(P_{1}(x_{1}))\cdot\cdot\cdot\mathrm{h}_{m}(P_{m}(x_{m}))
\]
by linearity of $\mathrm{h}_{m}.\;$We can write%
\[
\delta_{m}(\alpha)=\sum_{\sigma\in\mathcal{S}_{m}}(-1)^{l(\sigma)}%
x_{\sigma(1)}^{\alpha_{1}-\sigma(1)+1}(x_{\sigma(2)}^{\alpha_{2}-\sigma
(2)+2}+x_{\sigma(2)}^{\alpha_{2}-\sigma(2)})\cdot\cdot\cdot(x_{\sigma
(m)}^{\alpha_{m}-\sigma(m)+m}+x_{\sigma(m)}^{\alpha_{m}-\sigma(m)-m+2})
\]
and by the previous argument%
\[
\mathrm{h}_{m}(\delta_{m}(\alpha))=\sum_{\sigma\in\mathcal{S}_{m}%
}(-1)^{l(\sigma)}h_{\alpha_{1}-\sigma(1)+1}(h_{\alpha_{2}-\sigma
(2)+2}+h_{\alpha_{2}-\sigma(2)})\cdot\cdot\cdot(h_{\alpha_{m}-\sigma
(m)+m}+h_{\alpha_{m}-\sigma(m)-m+2})=u_{\alpha}%
\]
where the last equality follows from Proposition \ref{prop_fult}.\ By
(\ref{delda_factor}) we have $\delta_{m}(\alpha)=\mathrm{\omega}_{m}%
(x^{\alpha}).\;$Thus by applying $\mathrm{h}_{m}$ to this equality we obtain
$\mathrm{h}_{m}(\mathrm{\omega}_{m}(x^{\alpha}))=u_{\alpha}=\mathrm{u}%
_{m}(x^{\alpha})$.\ Hence $\mathrm{u}_{m}=\mathrm{h}_{m}\cdot\mathrm{\omega
}_{m}.\;$We derive the equality $\mathrm{u}_{m}=\mathrm{H}_{m}\cdot
\mathrm{\Omega}_{m}$ in a similar way starting from%
\[
\Delta_{m}(\alpha)=\sum_{\sigma\in\mathcal{S}_{m}}(-1)^{l(\sigma)}%
(x_{\sigma(1)}^{\alpha_{1}-\sigma(1)+1}+x_{\sigma(1)}^{\alpha_{2}-\sigma
(1)-1})\cdot\cdot\cdot(x_{\sigma(m)}^{\alpha_{m}-\sigma(m)+m}+x_{\sigma
(m)}^{\alpha_{m}-\sigma(m)-m}).
\]

$2:$ The arguments are the same than in $1$ once replacing the characters $h$
and $H$ respectively by the characters $e$ and $E.$
\end{proof}

\bigskip

Consider $\alpha=(\alpha_{1},...,\alpha_{m})\in\pi_{m}$ and two integers $i,j$
such that $1\leq i\leq j\leq m.$ The raising operator $R_{i,j}$ and the
lowering operator $L_{i,j}$ are respectively defined on $\pi_{m}$ by
$R_{i,j}(\alpha)=\alpha+\varepsilon_{i}-\varepsilon_{j}$ and $L_{i,j}%
(\alpha)=\alpha-\varepsilon_{i}-\varepsilon_{j}.$ From the previous lemma we obtain:

\begin{corollary}
\label{cor_s_h_R}For any partition $\mu=(\mu_{1},...,\mu_{m})$ we have
\begin{align*}
s_{\mu}  &  =\left(  \prod_{1\leq i<j\leq m}(1-R_{i,j})\prod_{1\leq r<s\leq
m}(1-L_{r,s})\right)  h_{\mu}\text{, }s_{\mu}=\left(  \prod_{1\leq i<j\leq
m}(1-R_{i,j})\prod_{1\leq r\leq s\leq m}(1-L_{r,s})\right)  H_{\mu},\\
s_{\mu}  &  =\left(  \prod_{1\leq i<j\leq n}(1-R_{i,j})\prod_{1\leq r<s\leq
n}(1-L_{r,s})\right)  e_{\mu^{\prime}},\text{ }s_{\mu}=\left(  \prod_{1\leq
i<j\leq n}(1-R_{i,j})\prod_{1\leq r\leq s\leq n}(1-L_{r,s})\right)
E_{\mu^{\prime}}%
\end{align*}
where $\mu^{\prime}=(\mu_{1}^{\prime},...,\mu_{n}^{\prime})$ is the conjugate
partition of $\mu.$
\end{corollary}

\begin{proof}
Let us write
\[
\phi_{m}=\prod_{1\leq i<j\leq m}(1-\frac{x_{i}}{x_{j}})\prod_{1\leq r<s\leq
m}(1-\frac{1}{x_{i}x_{j}})=\sum_{\alpha\in\pi_{m}}a(\alpha)x^{\alpha}.
\]
Then by $1$ of Proposition \ref{prop_s-h}, we have for any $\mu\in\pi_{m}%
^{+},$
\[
\mathrm{u}_{m}(x^{\mu})=\mathrm{h}_{m}\left(  \sum_{\alpha\in\pi_{m}}%
a(\alpha)x^{\alpha+\mu}\right)  =\sum_{\alpha\in\pi_{m}}a(\alpha)h_{\alpha
+\mu}=u_{\lambda}=s_{\lambda}%
\]
where the last equality follows from Proposition \ref{prop_fult}.\ This is
exactly equivalent to%
\[
s_{\mu}=\left(  \prod_{1\leq i<j\leq m}(1-R_{i,j})\prod_{1\leq r<s\leq
m}(1-L_{r,s})\right)  h_{\mu}.
\]
The arguments are essentially the same for the other equalities.
\end{proof}

\subsection{Expressions for the multiplicities of
representations\label{subsecC}}

Write
\[
\phi_{m}^{-1}=\sum_{\alpha\in\pi_{m}}f(\alpha)x^{\alpha}\text{ and }\Phi
_{m}^{-1}=\sum_{\alpha\in\pi_{m}}F(\alpha)x^{\alpha}.
\]
From Lemma \ref{prop_s-h} we deduce that $\mathrm{h}_{m}=\mathrm{u}_{m}%
\circ\mathrm{\omega}_{m}^{-1}$ and $\mathrm{H}_{m}=\mathrm{u}_{m}%
\circ\mathrm{\Omega}_{m}^{-1}.\;$By applying these identities to $x^{\mu}$
where $\mu$ is a partition of length $m$ with $\mu^{\prime}$ of length $n$ we
obtain as in Corollary \ref{cor_s_h_R}%
\begin{align*}
h_{\mu} &  =\left(  \prod_{1\leq i<j\leq m}\frac{1}{1-R_{i,j}}\prod_{1\leq
r<s\leq m}\frac{1}{1-L_{r,s}}\right)  s_{\mu}\text{, }H_{\mu}=\left(
\prod_{1\leq i<j\leq m}\frac{1}{1-R_{i,j}}\prod_{1\leq r\leq s\leq m}\frac
{1}{1-L_{r,s}}\right)  s_{\mu},\\
e_{\mu^{\prime}} &  =\left(  \prod_{1\leq i<j\leq n}\frac{1}{1-R_{i,j}}%
\prod_{1\leq r<s\leq n}\frac{1}{1-L_{r,s}}\right)  s_{\mu}\text{ and }%
E_{\mu^{\prime}}=\left(  \prod_{1\leq i<j\leq n}\frac{1}{1-R_{i,j}}%
\prod_{1\leq r\leq s\leq n}\frac{1}{1-L_{r,s}}\right)  s_{\mu}.
\end{align*}
These relations must be understood as a short way to write%
\begin{align*}
h_{\mu} &  =\sum_{\alpha\in\pi_{m}}f(\alpha)u_{\mu+\alpha}\text{, }H_{\mu
}=\sum_{\alpha\in\pi_{m}}F(\alpha)u_{\mu+\alpha},\\
e_{\mu^{\prime}} &  =\sum_{\beta\in\pi_{n}}f(\alpha)v_{\mu^{\prime}+\beta
}\text{ and }E_{\mu^{\prime}}=\sum_{\beta\in\pi_{n}}F(\alpha)v_{\mu^{\prime
}+\beta}.
\end{align*}
For any positive integer $l$ write $\rho_{l}=(l,l-1,...,1).$

\begin{proposition}
Consider a partition $\mu$ of length $m$ such that $\mu^{\prime}$ has length
$n.\;$Then for the three roots systems $B_{m},C_{m}$ and $D_{m}$ we have:%
\begin{align*}
\mathrm{(i)}  &  :\left\{
\begin{array}
[c]{c}%
h_{\mu}=\sum_{\lambda\in\pi_{m}}\sum_{\sigma\in\mathcal{S}_{m}}(-1)^{l(\sigma
)}f(\sigma(\lambda+\rho_{m})-\mu-\rho_{m})u_{\lambda}\\
H_{\mu}=\sum_{\lambda\in\pi_{m}}\sum_{\sigma\in\mathcal{S}_{m}}(-1)^{l(\sigma
)}F(\sigma(\lambda+\rho_{m})-\mu-\rho_{m})u_{\lambda}\text{ }%
\end{array}
\right.  ,\\
\mathrm{(ii)}  &  :\left\{
\begin{array}
[c]{c}%
e_{\mu^{\prime}}=\sum_{\nu\in\pi_{n}}\sum_{\sigma\in\mathcal{S}_{n}%
}(-1)^{l(\sigma)}f(\sigma(\nu+\rho_{n})-\mu^{\prime}-\rho_{n})v_{\nu}\\
E_{\mu^{\prime}}=\sum_{\nu\in\pi_{n}}\sum_{\sigma\in\mathcal{S}_{n}%
}(-1)^{l(\sigma)}F(\sigma(\nu+\rho_{n})-\mu^{\prime}-\rho_{n})v_{\nu}%
\end{array}
\right.  .
\end{align*}
\end{proposition}

\begin{proof}
$\mathrm{(i)}:$ Note first that the above relations do not depend on the root
system considered.\ Indeed for any nonnegative integer $m,$ we have
$\rho_{B_{n}}=\rho_{m}-(\frac{1}{2},....,\frac{1}{2}),\rho_{C_{m}}=\rho_{m}$
and $\rho_{D_{n}}=\rho_{m}-(1,...,1).\;$Thus $\sigma(\lambda+\rho_{B_{m}}%
)-\mu-\rho_{B_{m}}=\sigma(\lambda+\rho_{C_{m}})-\mu-\rho_{C_{m}}%
=\sigma(\lambda+\rho_{D_{m}})-\mu-\rho_{D_{m}}=\sigma(\lambda+\rho_{m}%
)-\mu-\rho_{m}.\;$We have%
\[
h_{\mu}=\sum_{\alpha\in\pi_{m}}f(\alpha)u_{\mu+\alpha}.
\]
From Lemma \ref{straight} we deduce that for any $\alpha\in\pi_{m}$ we have
$u_{\mu+\alpha}=0$ or there exits a partition $\lambda$ such that $\mu
+\alpha=\sigma(\lambda+\rho_{m})-\rho_{m}$ and $u_{\mu+\alpha}=(-1)^{l(\sigma
)}u_{\lambda}.\;$By setting $\alpha=\sigma(\lambda+\rho_{m})-\mu-\rho_{m}$ in
the above sum we obtain $h_{\mu}=\sum_{\lambda\in\pi_{m}}\sum_{\sigma
\in\mathcal{S}_{m}}(-1)^{l(\sigma)}f(\sigma(\lambda+\rho_{m})-\mu-\rho
_{m})u_{\lambda}.\;$The arguments are similar for the other assertions.
\end{proof}

\noindent From relations $\mathrm{(i)}$ and by using the fact that
$u_{\lambda}=s_{\lambda}$ for any partition $\lambda$ of length $m,$ we derive
the equalities%
\[
h_{\mu}=\sum_{\lambda\in\pi_{m}}u_{\lambda,\mu}s_{\lambda}\text{ and }H_{\mu
}=\sum_{\lambda\in\pi_{m}}U_{\lambda,\mu}s_{\lambda}%
\]
where%
\begin{equation}
u_{\lambda,\mu}=\sum_{\sigma\in\mathcal{S}_{m}}(-1)^{l(\sigma)}f(\sigma
(\lambda+\rho_{m})-\mu-\rho_{m})\text{ and }U_{\lambda,\mu}=\sum_{\sigma
\in\mathcal{S}_{m}}(-1)^{l(\sigma)}F(\sigma(\lambda+\rho_{m})-\mu-\rho
_{m})\label{U,u}%
\end{equation}
are respectively the multiplicities of $V(\lambda)$ in $\frak{h}(\mu)$ and
$\frak{H}(\mu).$ Note that $u_{\lambda,\mu}=0$ and $U_{\lambda,\mu}=0$ unless
$\left|  \mu\right|  \geq\left|  \lambda\right|  .$

\noindent For the relations $\mathrm{(ii)}$ the situation is more complicated
since the partitions $\nu$ obtained by applying straightening laws to the
$v_{\mu^{\prime}+\beta}$ yields polynomials $v_{\nu}$ where $\nu\in\pi_{n}%
^{+}$ is a partition of length $n$ so can not be necessarily regarded as the
conjugate partition of a partition $\lambda\in\pi_{m}^{+}.$ The straightening
law of Lemma \ref{straight} implies that $\left|  \nu\right|  =\left|
\mu^{\prime}\right|  .\;$Since $\left|  \mu\right|  =\left|  \mu^{\prime
}\right|  ,$ this problem disappear if we suppose $m\geq\left|  \mu\right|  $
since we will have $\nu_{1}\leq\left|  \nu\right|  \leq m$ and thus
$\nu^{\prime}\in\pi_{m}^{+}.\;$We can then set $\nu=\lambda^{\prime}$ with
$\lambda\in\pi_{m}$ and obtain%
\[
e_{\mu^{\prime}}=\sum_{\lambda\in\pi_{m}}v_{\lambda,\mu}s_{\lambda}\text{ and
}E_{\mu^{\prime}}=\sum_{\lambda\in\pi_{m}}V_{\lambda,\mu}s_{\lambda}.
\]
We deduce that
\begin{gather}
v_{\lambda,\mu}=u_{\lambda^{\prime},\mu^{\prime}}=\sum_{\sigma\in
\mathcal{S}_{n}}(-1)^{l(\sigma)}f(\sigma(\lambda^{\prime}+\rho_{n}%
)-\mu^{\prime}-\rho_{n})\label{v-V}\\
V_{\lambda,\mu}=U_{\lambda^{\prime},\mu^{\prime}}=\sum_{\sigma\in
\mathcal{S}_{n}}(-1)^{l(\sigma)}F(\sigma(\lambda^{\prime}+\rho_{n}%
)-\mu^{\prime}-\rho_{n})
\end{gather}
are respectively the multiplicities of $V(\lambda)$ in the tensor products
$\frak{e}(\mu)$ and $\frak{E}(\mu).$

\section{Quantification of the multiplicities}

\subsection{The functions $f_{q}$ and $F_{q}$}

Set%
\[
\phi_{m}(q)=\prod_{1\leq i<j\leq m}(1-q\frac{x_{i}}{x_{j}})\prod_{1\leq
r<s\leq m}(1-\frac{q}{x_{i}x_{j}})\text{ and }\Phi_{m}(q)=\prod_{1\leq i<j\leq
m}(1-q\frac{x_{i}}{x_{j}})\prod_{1\leq r\leq s\leq m}(1-\frac{q}{x_{i}x_{j}%
}).
\]
The functions $f_{q}$ and $F_{q}$ are obtained by considering the formal
series expansions of $\phi_{m}^{-1}(q)$ and $\Phi_{m}^{-1}(q).\;$Namely we
have
\begin{equation}
\phi_{m}^{-1}(q)=\sum_{\alpha\in\pi_{m}}f_{q}(\alpha)x^{\alpha}\text{ and
}\Phi_{m}^{-1}(q)=\sum_{\alpha\in\pi_{m}}F_{q}(\alpha)x^{\alpha}.
\label{def_f_F}%
\end{equation}

\subsection{Some $q$-analogues of multiplicities of $V(\lambda)$ in
$\frak{h}(\mu)$, $\frak{H}(\mu),$ $\frak{e}(\mu)$ or $\frak{E}(\mu)$}

Given $\lambda$ and $\mu$ two partitions of length $m,$ let $c_{\lambda,\mu
}(q)$ and $C_{\lambda,\mu}(q)$ be the two polynomials defined by%

\[
u_{\lambda,\mu}(q)=\sum_{\sigma\in\mathcal{S}_{m}}(-1)^{l(\sigma)}f_{q}%
(\sigma(\lambda+\rho_{m})-\mu-\rho_{m})\text{ and }U_{\lambda,\mu}%
(q)=\sum_{\sigma\in\mathcal{S}_{m}}(-1)^{l(\sigma)}F_{q}(\sigma(\lambda
+\rho_{m})-\mu-\rho_{m}).
\]
Then from the equalities (\ref{U,u}) and (\ref{v-V}) we obtain:

\begin{proposition}
\label{prop_multi}Let $\lambda$ and $\mu$ be two partitions of length $m.\;$Then

\begin{enumerate}
\item $u_{\lambda,\mu}(q)$ and $U_{\lambda,\mu}(q)$ are $q$-analogues of the
multiplicity of the representation $V(\lambda)$ in $\frak{h}(\mu)$ and
$\frak{H}(\mu),$

\item $v_{\lambda,\mu}(q)=u_{\lambda,^{\prime}\mu^{\prime}}(q)$ and
$V_{\lambda,\mu}(q)=U_{\lambda^{\prime},\mu^{\prime}}(q)$ are $q$-analogues of
the multiplicity of the representation $V(\lambda)$ in $\frak{e}(\mu)$ and
$\frak{E}(\mu)$ when the condition $m\geq\left|  \mu\right|  $ is satisfied.
\end{enumerate}
\end{proposition}

\noindent The following example is obtained from the explicit computation of
the function $f_{q}$ when $m=2.$

\begin{example}
Consider $\mu$ a partition of length $2$ and set $\mathcal{E}_{\mu}%
=\{\lambda\in\pi_{2}^{+},\lambda=(\mu_{1}+r-s,\mu_{2}-r-s),$ $s\in
\{0,...,\mu_{2}\},$ $r\in\{0,...,\mu_{2}-s\}\}.\;$Then for any partition
$\lambda$ of length $2$ we have:%
\[
u_{\lambda,\mu}(q)=\left\{
\begin{array}
[c]{c}%
q^{\mu_{1}-\lambda_{1}}\text{ if }\lambda\in\mathcal{E}_{\mu}\\
0\text{ otherwise}%
\end{array}
\right.  .
\]
\end{example}

\noindent\textbf{Remarks }

\noindent$\mathrm{(i):}$ It follows from the definition of the $q$-functions
$f_{q}$ and $F_{q}$ that $c_{\lambda,\mu}(q)=C_{\lambda,\mu}(q)=0$ if $\left|
\lambda\right|  >\left|  \mu\right|  .$

\noindent$\mathrm{(ii):}$ It is not trivial from the very definitions that
$u_{\lambda,\mu}(q)$ and $U_{\lambda,\mu}(q)$ are polynomials in $q$ with
nonnegative integer coefficients.\ This property will be proved in Section
\ref{sec_last} as a corollary of Theorem \ref{th_dual1}.

\section{\label{sec_last}The duality theorems}

\subsection{A duality theorem for the $q$-multiplicities in $\frak{h}(\mu)$
and $\frak{H}(\mu)$}

For any nonnegative integer $m,$ set $\kappa_{m}=(1,...,1)\in\pi_{m}.$

\begin{lemma}
\label{lem_sum_rest}Consider $\lambda,\mu$ two partitions of length $m$ such
that $\left|  \lambda\right|  \geq\left|  \mu\right|  .\;$Let $k$ be any
integer such that $k\geq\frac{\left|  \lambda\right|  -\left|  \mu\right|
}{2}$.\ Then we have
\begin{equation}
K_{\lambda+k\kappa_{m},\mu+k\kappa_{m}}(q)=\sum_{\sigma\in S_{m}%
}(-1)^{l(\sigma)}\mathcal{P}_{q}(\sigma(\lambda+\rho_{m})-(\mu+\rho_{m}))
\label{sum-rest}%
\end{equation}
where the sum is indexed by the elements of the symmetric group $S_{m}.$
\end{lemma}

\begin{proof}
Since $\mathcal{P}_{q}(\alpha)=0$ if $\alpha$ is not a linear combination of
positive roots with nonnegative coefficients, we have $\mathcal{P}_{q}%
(\alpha)=0$ for any $\alpha\in\pi_{m}$ such that $\left|  \alpha\right|
<0.\;$Consider $\delta=(\delta_{1},...,\delta_{m})\in\pi_{m}$ and $w\in
W_{n}.\;$Write $w(\delta)=(\delta_{1}^{w},...,\delta_{m}^{w})$ and denote by
$E_{w,\delta}=\{i_{1},...,i_{p}\}$ the set of the indices $i_{p}$ such that
$\delta_{i_{k}}$ and $\delta_{i_{k}}^{w}$ have opposite signs.\ Define the sum
$S_{w,\delta}=\sum_{i_{p}\in E_{w,\delta}}\delta_{i_{k}}.\;$Then $\left|
w(\delta)\right|  =\left|  \delta\right|  -2S_{w,\delta}$.\ Now consider $k$ a
nonnegative integer and set $\delta=(\lambda+\rho_{m}+k\kappa_{m}).\;$We have
$\left|  w(\lambda+\rho_{m}+k\kappa_{m})\right|  =\left|  (\lambda+\rho
_{m}+k\kappa_{m})\right|  -2S_{w,\delta}.\;$But $S_{w,\delta}=S_{w,\lambda
+\rho_{m}}+kp.\;$Thus we obtain%
\begin{multline*}
\left|  w(\lambda+\rho_{m}+k\kappa_{m})-(\mu+\rho_{m}+k\kappa_{m})\right|
=\left|  (\lambda+\rho_{m}+k\kappa_{m})\right|  -2S_{w,\lambda+\rho_{m}%
}-\left|  (\mu+\rho_{m}+k\kappa_{m})\right|  -2kp=\\
\left|  \lambda\right|  -\left|  \mu\right|  -2S_{w,\lambda+\rho_{m}}-2kp.
\end{multline*}
When $w\notin\mathcal{S}_{m},$ we have $p\geq1\;$and $S_{w,\lambda+\rho_{m}%
}\geq1$ since the coordinates of $\lambda+\rho_{m}$ are all positive. Hence
$\left|  w(\lambda+\rho_{m}+k\kappa_{m})-(\mu+\rho_{m}+k\kappa_{m})\right|
<\left|  \lambda\right|  -\left|  \mu\right|  -2k$ and is negative as soon as
$k\geq\frac{\left|  \lambda\right|  -\left|  \mu\right|  }{2}.$ For such an
integer $k$ the sum defining $K_{\lambda+k\kappa_{m},\mu+k\kappa_{m}}(q)$
normally running on $W_{m}$ can be restricted to (\ref{sum-rest}) and we
obtain%
\[
K_{\lambda+k\kappa_{m},\mu+k\kappa_{m}}(q)=\sum_{\sigma\in S_{m}%
}(-1)^{l(\sigma)}\mathcal{P}_{q}(\sigma(\lambda+\rho_{m}+k\kappa_{m}%
)-(\mu+\rho_{m}+k\kappa_{m})).
\]
Since $\sigma\in\mathcal{S}_{m},$ we have $\sigma(k\kappa_{m})=k\kappa_{m}.$
Thus
\[
K_{\lambda+k\kappa_{m},\mu+k\kappa_{m}}(q)=\sum_{\sigma\in S_{m}%
}(-1)^{l(\sigma)}\mathcal{P}_{q}(\sigma(\lambda+\rho_{m})-(\mu+\rho_{m})).
\]
\end{proof}

\bigskip

We define the involution $I$ on $\pi_{m}$ by $I(\alpha_{1},...,\alpha
_{m})=(-\alpha_{m},...,-\alpha_{1})$ for any $\alpha=(\alpha_{1}%
,...,\alpha_{m})\in\pi_{m}.$

\begin{lemma}
\label{lem_invol_I}For any $\alpha=(\alpha_{1},...,\alpha_{m})\in\pi_{m}$ we
have%
\[
f_{q}(\alpha)=\mathcal{P}_{q}^{D_{m}}(I(\alpha))\text{ and }F_{q}%
(\alpha)=\mathcal{P}_{q}^{C_{m}}(I(\alpha))
\]
where $\mathcal{P}_{q}^{B_{m}}$ and $\mathcal{P}_{q}^{D_{m}}$ are the
$q$-Kostant's partition functions associated respectively to the root systems
$B_{m}$ and $D_{m}.$
\end{lemma}

\begin{proof}
By abuse of notation we also denote by $I$ the ring automorphism of
$\mathcal{L}_{m}$ defined by $I(x^{\alpha})=x^{I(\alpha)}.$ The image of the
root systems $C_{m}$ and $D_{m}$ by $I$ are respectively
\begin{equation}
\left\{
\begin{tabular}
[c]{l}%
$\{\varepsilon_{i}-\varepsilon_{j},-\varepsilon_{i}-\varepsilon_{j}\text{ with
}1\leq i<j\leq m\}\cup\{-2\varepsilon_{i}\text{ with }1\leq i\leq m\}\text{
for the root system }C_{m}$\\
$\{\varepsilon_{i}-\varepsilon_{j},-\varepsilon_{i}-\varepsilon_{j}\text{ with
}1\leq i<j\leq m\}\text{ for the root system }D_{m}$%
\end{tabular}
\right.  .\label{Iroot}%
\end{equation}
By applying $I$ to the equality%
\[
\prod_{\alpha\in R_{C_{m}}^{+}}\dfrac{1}{1-qx^{\alpha}}=\sum_{\beta\in\pi_{m}%
}\mathcal{P}_{q}^{C_{m}}(\beta)x^{\beta}%
\]
we obtain
\[
\prod_{1\leq i<j\leq m}\frac{1}{(1-q\frac{x_{i}}{x_{j}})}\prod_{1\leq r\leq
s\leq m}\frac{1}{(1-\frac{q}{x_{r}x_{s}})}=\sum_{\beta\in\pi_{m}}%
\mathcal{P}_{q}^{C_{m}}(\beta)x^{I(\beta)}.
\]
Set $\alpha=I(\beta).\;$The equality becomes%
\[
\Phi_{m}^{-1}(q)=\sum_{\alpha\in\pi_{m}}\mathcal{P}_{q}^{C_{m}}(I(\alpha
))x^{\alpha}%
\]
and from the definition (see \ref{def_f_F}) of the function $F_{q},$ we obtain
$\mathcal{P}_{q}^{B_{m}}(I(\alpha))=F_{q}(\alpha).$ The assertion with $f_{q}$
is proved in the same way by considering the root system $D_{m}.$
\end{proof}

\bigskip

\noindent Given $\sigma\in\mathcal{S}_{m},$ denote by $\sigma^{\ast}$ the
permutation defined by%
\[
\sigma^{\ast}(k)=\sigma(m-k+1).
\]
For any $i\in\{1,...,m-1\},$ we have $s_{i}^{\ast}=s_{m-i}$ so that
$l(\sigma)=l(\sigma^{\ast}).\;$The following Lemma is straightforward:

\begin{lemma}
\label{lemma-invol_Sn}The map $\sigma\rightarrow\sigma^{\ast}$ is an
involution of the group $\mathcal{S}_{m}$ and we have $\sigma(I(\beta
))=I(\sigma^{\ast}(\beta))$ for any $\beta\in\pi_{m},\sigma\in\mathcal{S}_{m}.$
\end{lemma}

\begin{lemma}
Let $\lambda,\mu$ two partitions of length $m$ and $\sigma\in\mathcal{S}_{m}.$
Then%
\begin{gather*}
(-1)^{l(\sigma)}f_{q}(\sigma(\lambda+\rho_{m})-(\mu+\rho_{m}))=(-1)^{l(\sigma
^{\ast})}\mathcal{P}_{q}^{D_{m}}(\sigma^{\ast}(I(\lambda)+\rho_{m}%
)-(I(\mu)+\rho_{m}))\text{ and}\\
(-1)^{l(\sigma)}F_{q}(\sigma(\lambda+\rho_{m})-(\mu+\rho))=(-1)^{l(\sigma
^{\ast})}\mathcal{P}_{q}^{C_{m}}(\sigma^{\ast}(I(\lambda)+\rho_{m}%
)-(I(\mu)+\rho_{m})).
\end{gather*}
\end{lemma}

\begin{proof}
Since $l(\sigma)=l(\sigma^{\ast}),$ it suffices to prove the equalities
\begin{gather*}
f_{q}(\sigma(\lambda+\rho_{m})-(\mu+\rho_{m}))=\mathcal{P}_{q}^{D_{m}}%
(\sigma(I(\lambda)+\rho_{m})-(I(\lambda)+\rho_{m}))\text{ and}\\
F_{q}(\sigma(\lambda+\rho_{m})-(\mu+\rho_{m}))=\mathcal{P}_{q}^{C_{m}}%
(\sigma(I(\mu)+\rho_{m})-(I(\mu)+\rho_{m})).
\end{gather*}
Set $P=\mathcal{P}_{q}^{C_{m}}(\sigma^{\ast}(I(\lambda)+\rho_{m})-(I(\mu
)+\rho_{m})).\;$From the above Lemma we deduce%
\[
P=\mathcal{P}_{q}^{C_{m}}(I(\sigma(\lambda)+\sigma^{\ast}(\rho_{m}%
)-I(\mu)-\rho_{m}).
\]
Now an immediate computation shows that $\sigma^{\ast}(\rho_{m})-\rho
_{m}=I(\sigma(\rho_{m})-\rho_{m}).\;$Thus we derive%
\[
P=\mathcal{P}_{q}^{C_{m}}(I(\sigma(\lambda+\rho_{m})-\mu-\rho_{m}%
))=F_{q}(\sigma(\lambda+\rho_{m})-\mu-\rho_{m})
\]
where the last equality follows from Lemma \ref{lem_invol_I}.

\noindent We obtain the equality $f_{q}(\sigma(\lambda+\rho_{m})-(\mu+\rho
_{m}))=\mathcal{P}_{q}^{D_{m}}(\sigma(I(\lambda)+\rho_{m})-(I(\lambda
)+\rho_{m}))$ in a similar way.
\end{proof}

\begin{theorem}
\label{th_dual1}Consider $\lambda,\mu$ two partitions of length $m$ and set
$n=\max(\lambda_{1},\mu_{1})$.\ Let $k$ be any integer such that $k\geq
\frac{\left|  \mu\right|  -\left|  \lambda\right|  }{2}.\;$Then $\widehat
{\lambda}=(n-\lambda_{m},...,n-\lambda_{1})$ and $\widehat{\mu}=(n-\mu
_{m},...,n-\mu_{1})$ are partitions of length $m$ and
\[
\left\{
\begin{tabular}
[c]{l}%
$u_{\lambda,\mu}(q)=K_{\widehat{\lambda}+k\kappa_{m},\widehat{\mu}+k\kappa
_{m}}^{D_{m}}(q)$\\
$U_{\lambda,\mu}(q)=K_{\widehat{\lambda}+k\kappa_{m},\widehat{\mu}+k\kappa
_{m}}^{C_{m}}(q)$%
\end{tabular}
\right.
\]
\end{theorem}

\begin{proof}
First $\widehat{\lambda}$ and $\widehat{\mu}$ are clearly partitions of length
$m$ since $n=\mathrm{max}(\lambda_{1},\mu_{1}).$ It follows from the
definition of $U_{\lambda,\mu}(q)$ and the above lemma that%
\[
U_{\lambda,\mu}(q)=\sum_{\sigma\in\mathcal{S}_{m}}(-1)^{l(\sigma)}F_{q}%
(\sigma(\lambda+\rho_{m})-\mu-\rho_{m})=\sum_{\sigma^{\ast}\in\mathcal{S}_{m}%
}(-1)^{l(\sigma^{\ast})}\mathcal{P}_{q}^{C_{m}}(\sigma^{\ast}(I(\lambda
)+\rho_{m}))-(I(\mu)+\rho_{m})).
\]
Then by Lemma \ref{lemma-invol_Sn} we obtain%
\[
U_{\lambda,\mu}(q)=\sum_{\sigma\in\mathcal{S}_{m}}(-1)^{l(\sigma)}%
\mathcal{P}_{q}^{C_{m}}(\sigma(I(\lambda)+\rho_{m}))-(I(\mu)+\rho_{m})).
\]
We have $\sigma(\lambda^{\ast}+\rho_{m}+n\kappa_{m})=\sigma(\lambda^{\ast
}+\rho_{m})+n\kappa_{m}$ since $\sigma\in\mathcal{S}_{m}.$ So we can write%
\[
U_{\lambda,\mu}(q)=\sum_{\sigma\in\mathcal{S}_{m}}(-1)^{l(\sigma)}%
\mathcal{P}_{q}^{C_{m}}(\sigma(I(\lambda)+n\kappa_{m}+\rho_{m}))-(I(\mu
)+n\kappa_{m}+\rho_{m})).
\]
Since $\widehat{\lambda}=I(\lambda)+n\kappa_{m}$ and $\widehat{\mu}%
=I(\mu)+n\kappa_{m}$ we derive%
\[
U_{\lambda,\mu}(q)=\sum_{\sigma\in\mathcal{S}_{m}}(-1)^{l(\sigma)}%
\mathcal{P}_{q}^{C_{m}}(\sigma(\widehat{\lambda}+\rho_{m})-(\widehat{\mu}%
+\rho_{m}))=K_{\widehat{\lambda}+k\kappa_{m},\widehat{\mu}+k\kappa_{m}}%
^{C_{m}}(q)
\]
by Lemma \ref{lem_sum_rest}.

\noindent We obtain similarly the equality $u_{\lambda,\mu}(q)=K_{\widehat
{\lambda}+k\kappa_{m},\widehat{\mu}+k\kappa_{m}}^{D_{m}}(q)$ by replacing
$\mathcal{P}_{q}^{C_{m}}$ by $\mathcal{P}_{q}^{D_{m}}.$
\end{proof}

\begin{example}
Consider $\mu=(4,2,1)$ and $\lambda=(2,1,0).\;$We have $n=4$, $\widehat{\mu
}=(3,2,0)$ and $\widehat{\lambda}=(4,3,2).\;$We choose $k=2.\;$Then we obtain
the equalities%
\[
\left\{
\begin{tabular}
[c]{l}%
$u_{\lambda,\mu}(q)=K_{(6,5,4),(5,4,2)}^{D_{m}}(q)=q^{3}+q^{2}$\\
$U_{\lambda,\mu}(q)=K_{(6,5,4),(5,4,2))}^{C_{m}}(q)=q^{5}+2q^{4}+3q^{3}%
+2q^{2}$%
\end{tabular}
\right.  .
\]
\end{example}

\noindent By using the fact that the Kostka-Foulkes polynomials have
nonnegative integer coefficients \cite{Lu} we obtain the following corollary.

\begin{corollary}
The polynomials $u_{\lambda,\mu}(q)$ and $U_{\lambda,\mu}(q)$ have nonnegative
integers coefficients.
\end{corollary}

\noindent We also recover a property of the Kostka-Foulkes polynomials
associated to the root system $A_{m}$ proved in \cite{LSc1}.

\begin{corollary}
\label{coe-dual-A}Consider $\lambda,\mu$ two partitions of length $m$ such
that $\left|  \lambda\right|  =\left|  \mu\right|  $ and set $n=\max
(\lambda_{1},\mu_{1})$.\ Then the Kostka-Foulkes polynomials associated to the
root system $A_{m}$ verifies%
\[
K_{\lambda,\mu}^{A_{m}}(q)=K_{\widehat{\lambda},\widehat{\mu}}^{A_{m}}(q)
\]
where $\widehat{\lambda}=(n-\lambda_{m},...,n-\lambda_{1})$ and $\widehat{\mu
}=(n-\mu_{m},...,n-\mu_{1}).$
\end{corollary}

\begin{proof}
Suppose that $\beta$ is a linear combination of $I(R_{C_{m}}^{+})$ with
nonnegative coefficients such that $\left|  \beta\right|  =0.\;$Then $\beta$
is necessarily a linear combination of the roots $\varepsilon_{i}%
-\varepsilon_{j},1\leq$ $i<j\leq m$ with nonnegative coefficients (see
(\ref{Iroot})) that is, a linear combination with nonnegative coefficients of
the positive roots associated to the root system $A_{m}$.\ This implies that
\[
f_{q}(\beta)=F_{q}(\beta)=\mathcal{P}_{q}^{A_{m}}(\beta)
\]
where $\mathcal{P}_{q}^{A_{m}}$ is the $q$-Kostant's partition function
associated to the root system $A_{m}.\;$For any $\sigma\in\mathcal{S}_{m},$ we
have $\left|  \sigma(\lambda+\rho_{m})-(\mu+\rho_{m})\right|  =0$ since
$\left|  \lambda\right|  =\left|  \mu\right|  .\;$Thus
\[
f_{q}(\sigma(\lambda+\rho_{m})-(\mu+\rho_{m}))=F_{q}(\sigma(\lambda+\rho
_{m})-(\mu+\rho_{m}))=\mathcal{P}_{q}^{A_{m}}(\sigma(\lambda+\rho_{m}%
)-(\mu+\rho_{m}))
\]
and the multiplicities $u_{\lambda,\mu}(q)$ and $U_{\lambda,\mu}(q)$ coincide
with the Kostka-Foulkes polynomial $K_{\lambda,\mu}^{A_{m}}(q)$ when $\left|
\lambda\right|  =\left|  \mu\right|  .\;$Moreover by applying Theorem
\ref{th_dual1} with $\left|  \lambda\right|  =\left|  \mu\right|  $ and $k=0,$
we obtain $U_{\lambda,\mu}(q)=K_{\widehat{\lambda},\widehat{\mu}}^{C_{m}%
}(q)=K_{\widehat{\lambda},\widehat{\mu}}^{A_{m}}(q)$ where the last equality
is due to the fact that the Kostka-Foulkes polynomials of types $B_{m},C_{m}$
or $D_{m}$ are Kostka-Foulkes polynomials associated to the root system
$A_{m}$ when $\left|  \lambda\right|  =\left|  \mu\right|  .\;$So we derive
the equality $K_{\lambda,\mu}^{A_{m}}(q)=K_{\widehat{\lambda},\widehat{\mu}%
}^{A_{m}}(q).$
\end{proof}

We have seen that $U_{\lambda,\mu}(q)$ can be regarded as a $q$-analogue of
the multiplicity of the representation $V(\lambda)$ in$\frak{\;H}^{C_{m}}%
(\mu).\;$In \cite{Ok}, Hatayama, Kuniba, Okado and Takagi have introduced
another quantification $X_{\lambda,\mu}(q)$ of this multiplicity based on the
determination of the combinatorial $R$ matrix of the $U_{q}^{\prime}%
(C_{m}^{(1)})$-crystals $B_{k}$.\ Considered as the crystal graph of the
$U_{q}(C_{m})$-module $M_{k},$ $B_{k}$ can be identify with
\[
B(k\Lambda_{1})\oplus B((k-2)\Lambda_{1}^{C_{m}})\oplus\cdot\cdot\cdot
\cdot\oplus B(k\operatorname{mod}2\Lambda_{1}^{C_{m}})
\]
where for any $i\in\{k,k-2,...,k\mathrm{mod}2\}$, $B(\ k\Lambda_{1}^{C_{m}})$
is the graph corresponding to the irreducible finite dimensional highest
weight $U_{q}(C_{m})$-module of highest weight $k\Lambda_{1}^{C_{m}}.\;$Note
that the character of $M_{k}$ is equal to $H_{k}^{C_{m}}.$

\noindent Recall that the combinatorial $R$-matrix associated to crystals
$B_{k}$ is equivalent to the description of the crystal graph isomorphisms
\[
\left\{
\begin{array}
[c]{c}%
B_{l}\otimes B_{k}\overset{\simeq}{\rightarrow}B_{k}\otimes B_{l}\\
b_{1}\otimes b_{2}\longmapsto b_{2}^{\prime}\otimes b_{1}^{\prime}%
\end{array}
\right.
\]
together with the energy function $H$ on $B_{l}\otimes B_{k}.\;$The
multiplicity of $V(\lambda)$ in$\frak{\;H}^{C_{m}}(\mu)$ is then equal to the
number of highest weight vertices of weight $\lambda$ in the crystal $B_{\mu
}=B_{\mu_{1}}\otimes\cdot\cdot\cdot\otimes B_{\mu_{m}}$.\ Then $X_{\lambda
,\mu}(q)$ is defined by
\[
X_{\lambda,\mu}(q)=\sum_{b\in E_{\lambda}}q^{\sum_{0\leq i<j\leq m}%
H(b_{i}\otimes b_{j}^{(i+1)})}%
\]
where $E_{\lambda}$ is the set of highest weight vertices $b=b_{1}\otimes
\cdot\cdot\cdot\otimes b_{m}$ in $B_{\mu}$ of highest weight $\lambda,$
$b_{j}^{(i)}$ is determined by the crystal isomorphism%
\[%
\begin{array}
[c]{c}%
B_{\mu_{i}}\otimes B_{\mu_{i+1}}\otimes B_{\mu_{i+2}}\otimes\cdot\cdot
\cdot\otimes B_{\mu_{j}}\rightarrow B_{\mu_{i}}\otimes B_{\mu_{j}}\otimes
B_{\mu_{i+1}}\cdot\cdot\cdot\otimes B_{\mu_{j-1}}\\
b_{i}\otimes b_{i+1}\otimes\cdot\cdot\cdot\otimes b_{j}\rightarrow b_{j}%
^{(i)}\otimes b_{i}^{\prime}\otimes\cdot\cdot\cdot\otimes b_{j-1}^{\prime}%
\end{array}
\]
and for any $j=1,...,m,$ $H(b_{0}\otimes b_{j}^{(1)})$ depends only on
$b_{j}^{(1)}$.

\noindent Many computations suggest the following conjecture

\begin{conjecture}
For any partition $\lambda$ and $\mu$ of length $m$ with $\left|  \mu\right|
\geq\left|  \lambda\right|  $%
\[
U_{\lambda,\mu}(q)=q^{\left|  \mu\right|  -\left|  \lambda\right|  }%
X_{\lambda,\mu}(q).
\]
\end{conjecture}

\noindent Note that the conjecture is in particular true for all the examples
given in the tables of \cite{Ok}.\ 

\subsection{A duality theorem for the $q$-multiplicities in $\frak{e}(\mu)$
and $\frak{E}(\mu)$}

Consider $\lambda,\mu$ two partitions of length $l$ such that $l\geq\left|
\mu\right|  \geq\left|  \lambda\right|  .\;$Write $m=\mathrm{max}(\lambda
_{1},\mu_{1})$.\ Then by adding to $\lambda^{\prime}$ and $\mu^{\prime}$ the
required numbers of parts $0$ we can consider them as partitions of length
$m.\;$Set $n=\mathrm{max}(\lambda_{1}^{\prime},\mu_{1}^{\prime}).\;$We define
the partitions $\widetilde{\lambda}$ and $\widetilde{\mu}$ belonging to
$\pi_{m}$ by $\widetilde{\lambda}=(n-\lambda_{m}^{\prime},...,n-\lambda
_{1}^{\prime})$ and $\widetilde{\mu}=(n-\mu_{m}^{\prime},...,n-\mu_{1}%
^{\prime})$.

\begin{theorem}
\label{th_dual2}With the above notations, we have for any integer $k\geq
\frac{\left|  \mu\right|  -\left|  \lambda\right|  }{2}$%
\[
\left\{
\begin{tabular}
[c]{l}%
$\mathrm{(i)}:v_{\lambda,\mu}(q)=K_{\widetilde{\lambda}+k\kappa_{m}%
,\widetilde{\mu}+k\kappa_{m}}^{D_{m}}(q)$\\
$\mathrm{(ii)}:V_{\lambda,\mu}(q)=K_{\widetilde{\lambda}+k\kappa
_{m},\widetilde{\mu}+k\kappa_{m}}^{C_{m}}(q)$%
\end{tabular}
\right.  .
\]
\end{theorem}

\begin{proof}
Since $l\geq\left|  \mu\right|  ,$ we have by Proposition \ref{prop_multi} the
equality $v_{\lambda,\mu}(q)=u_{\lambda^{\prime},\mu^{\prime}}(q).$ Moreover
we have $n\geq\max(\lambda_{1}^{\prime},\mu_{1}^{\prime})$ and $k\geq
\frac{\left|  \mu^{\prime}\right|  -\left|  \lambda^{\prime}\right|  }{2}$ for
$\left|  \lambda^{\prime}\right|  =\left|  \lambda\right|  $ and $\left|
\mu^{\prime}\right|  =\left|  \mu\right|  .\;$Hence by applying Theorem
\ref{th_dual1}\ we obtain $v_{\lambda,\mu}(q)=K_{\widehat{\lambda^{\prime}%
}+k\kappa_{m},\widehat{\mu^{\prime}}+k\kappa_{m}}^{D_{m}}(q)$ where
$\widehat{\lambda^{\prime}}=(n-\lambda_{m}^{\prime},...,n-\lambda_{1}^{\prime
})=\widetilde{\lambda}$ and $\widehat{\mu^{\prime}}=(n-\mu_{m}^{\prime
},...,n-\mu_{1}^{\prime})=\widetilde{\mu}.$ So $\mathrm{(i)}$ is proved$.$ We
obtain $\mathrm{(ii)}$ similarly.
\end{proof}

\begin{example}
For $\lambda=(2,1,0,0,0)$ and $\mu=(2,2,1,0,0)$ we have $l=5,m=2$.\ Moreover
$\lambda^{\prime}=(2,1)$, $\mu^{\prime}=(3,2)$ and $n=3.\;$So $\widetilde
{\lambda}=(2,1)$ and $\widetilde{\mu}=(1,0)$.\ Hence for $k=1$%
\[
\left\{
\begin{tabular}
[c]{l}%
$\mathrm{(i)}:v_{\lambda,\mu}(q)=K_{(3,2),(2,1)}^{D_{m}}(q)=q$\\
$\mathrm{(ii)}:V_{\lambda,\mu}(q)=K_{(3,2),(2,1)}^{C_{m}}(q)=q^{2}+q$%
\end{tabular}
\right.
\]
\end{example}

\noindent\textbf{Remark }When $\lambda,\mu$ are considered as weight
associated to the root system $C_{l},$ the above theorem is essentially the
quantification of a duality result explicited by Foulle \cite{Fou} from
results of \cite{HO} for the dual pair $(Sp(2l),Sp(2m)).$

\section{Multiplicities for types $B_{m},C_{m},D_{m}$ and Kostka numbers}

\subsection{A relations between $q$-Kostant's partition functions}

Consider a nonnegative integer $k$ and define the finite sets%
\[
\left\{
\begin{tabular}
[c]{l}%
$\mathcal{C}_{k}^{m}=\{\beta\in\pi_{m},\beta=\sum_{1\leq r\leq s\leq m}%
e_{r,s}(\varepsilon_{r}+\varepsilon_{s})$ with $e_{r,s}\geq0$ and $\left|
\beta\right|  =2k\}$\\
$\mathcal{D}_{k}^{m}=\{\beta\in\pi_{m},\beta=\sum_{1\leq r<s\leq m}%
e_{r,s}(\varepsilon_{r}+\varepsilon_{s})$ with $e_{r,s}\geq0$ and $\left|
\beta\right|  =2k\}$%
\end{tabular}
\right.  .
\]
Note that each $\beta\in\mathcal{C}_{k}^{m}$ (resp. $\beta\in\mathcal{D}%
_{k}^{m})$ verifies $\left|  \beta\right|  =2\sum_{1\leq r\leq s\leq m}%
e_{r,s}$ (resp.\ $\left|  \beta\right|  =2\sum_{1\leq r<s\leq m}e_{r,s}%
).\;$This implies that%

\[
\prod_{1\leq r\leq s\leq m}\frac{1}{(1-qx_{r}x_{s})}=\sum_{k\geq0}\sum
_{\beta\in\mathcal{C}_{k}^{m}}c_{\beta}^{C_{m}}q^{k}x^{\beta}\text{
and\ }\prod_{1\leq r<s\leq m}\frac{1}{(1-qx_{r}x_{s})}=\sum_{k\geq0}%
\sum_{\beta\in\mathcal{C}_{k}^{m}}c_{\beta}^{D_{m}}q^{k}x^{\beta}%
\]
where $c_{\beta}^{C_{m}}$ (resp.\ $c_{\beta}^{D_{m}})$ is the number of ways
to decompose $\beta$ as $\beta=\sum_{1\leq r\leq s\leq m}e_{r,s}%
(\varepsilon_{r}+\varepsilon_{s})$ (resp.\ $\beta=\sum_{1\leq r<s\leq
m}e_{r,s}(\varepsilon_{r}+\varepsilon_{s}))$ with $e_{r,s}\geq0$.

\begin{lemma}
For any $\beta\in\pi_{m}$ with $\left|  \beta\right|  =2k\geq0,$ we have%
\[
\mathcal{P}_{q}^{C_{m}}(\beta)=\sum_{\delta\in\mathcal{C}_{k}^{m}}c_{\delta
}^{C_{m}}q^{k}\mathcal{P}_{q}^{A_{m}}(\beta-\delta)\text{ and }\mathcal{P}%
_{q}^{D_{m}}(\beta)=\sum_{\delta\in\mathcal{D}_{k}^{m}}c_{\delta}^{D_{m}}%
q^{k}\mathcal{P}_{q}^{A_{m}}(\beta-\delta).
\]
\end{lemma}

\begin{proof}
We have:%
\[
\prod_{1\leq i<j\leq m}\frac{1}{(1-q\frac{x_{i}}{x_{j}})}\prod_{1\leq r\leq
s\leq m}\frac{1}{(1-\frac{q}{x_{r}x_{s}})}=\sum_{\eta\in\pi_{m}}\sum
_{\delta\in\pi_{m}}c_{\delta}^{C_{m}}q^{\left|  \delta\right|  /2}%
\mathcal{P}_{q}^{A_{m}}(\eta)x^{\delta+\eta}%
\]
which implies the equality $\mathcal{P}_{q}^{C_{m}}(\beta)=\sum_{\eta
+\delta=\beta}c_{\delta}^{C_{m}}q^{\left|  \delta\right|  /2}\mathcal{P}%
_{q}^{A_{m}}(\eta).\;$Since $\mathcal{P}_{q}^{A_{m}}(\eta)=0$ when $\left|
\eta\right|  \neq0,$ we can suppose $\left|  \eta\right|  =0$ and $\left|
\delta\right|  =\left|  \beta\right|  $ in the previous sum.\ Then $\delta
\in\mathcal{C}_{k}^{m}$ and the result follows immediately.\ The proof for
$\mathcal{P}_{q}^{D_{m}}(\beta)$ is similar.
\end{proof}

\subsection{Multiplicities in terms of Kostka numbers}

\noindent Suppose that $\lambda$ is a partition of length $m$ and consider
$\gamma$ belonging to $\pi_{m}.\;$Then we can define the polynomial%
\[
K_{\lambda,\gamma}^{A_{m}}(q)=\sum_{\sigma\in\mathcal{S}_{m}}(-1)^{l(\sigma
)}\mathcal{P}_{q}^{A_{m}}(\sigma(\lambda+\rho_{m})-(\gamma+\rho_{m})).
\]
Note that the coefficients of $K_{\lambda,\gamma}^{A_{m}}(q)$ may be negative
when $\gamma$ is not a partition.\ Nevertheless $K_{\lambda,\gamma}^{A_{m}%
}=K_{\lambda,\gamma}^{A_{m}}(1)$ is equal to the dimension of the weight space
of weight $\gamma$ in $V^{A_{m}}(\lambda).$

\begin{proposition}
Consider $\lambda,\mu$ two partitions of length $m$ such that $k=\left|
\mu\right|  -\left|  \lambda\right|  \geq0.\;$Define $\widehat{\lambda}$ and
$\widehat{\mu}$ as in Theorem \ref{th_dual1}.\ Then%
\[
u_{\lambda,\mu}(q)=\sum_{\delta\in\mathcal{D}_{k}^{m}}c_{\delta}^{D_{m}%
}q^{\frac{\left|  \mu\right|  -\left|  \lambda\right|  }{2}}K_{\widehat
{\lambda},\widehat{\mu}+\delta}^{A_{m}}(q)\text{ and }U_{\lambda,\mu}%
(q)=\sum_{\delta\in\mathcal{C}_{k}^{m}}c_{\delta}^{C_{m}}q^{\frac{\left|
\mu\right|  -\left|  \lambda\right|  }{2}}K_{\widehat{\lambda},\widehat{\mu
}+\delta}^{A_{m}}(q)\text{.}%
\]
\end{proposition}

\begin{proof}
We have seen in the proof of Theorem \ref{th_dual1} that%
\[
U_{\lambda,\mu}=\sum_{\sigma\in\mathcal{S}_{m}}(-1)^{l(\sigma)}\mathcal{P}%
_{q}^{C_{m}}(\sigma(\widehat{\lambda}+\rho_{m})-(\widehat{\mu}+\rho_{m})).
\]
Hence from the above lemma we derive%
\[
U_{\lambda,\mu}=\sum_{\delta\in\mathcal{C}_{k}^{m}}c_{\delta}^{C_{m}%
}q^{\left|  \delta\right|  /2}\sum_{\sigma\in\mathcal{S}_{m}}(-1)^{l(\sigma
)}\mathcal{P}_{q}^{A_{m}}(\sigma(\widehat{\lambda}+\rho_{m})-(\widehat{\mu
}+\delta+\rho_{m}))
\]
which yields the desired result since $K_{\widehat{\lambda},\widehat{\mu
}+\delta}^{A_{m}}(q)=\sum_{\sigma\in\mathcal{S}_{m}}(-1)^{l(\sigma
)}\mathcal{P}_{q}^{A_{m}}(\sigma(\widehat{\lambda}+\rho_{m})-(\widehat{\mu
}+\delta+\rho_{m})).$
\end{proof}

\noindent From the above proposition and Proposition \ref{prop_multi} we
deduce by setting $q=1:$

\begin{corollary}
With the notations of Theorems \ref{th_dual1} and \ref{th_dual2} the following
assertions holds.

\begin{itemize}
\item  Consider $\lambda,\mu$ two partitions of length $m$ such that
$k=\frac{\left|  \mu\right|  -\left|  \lambda\right|  }{2}\in\mathbb{N}%
$.\ Then for types $B_{m},C_{m},D_{m},$ the multiplicity of the representation
$V(\lambda)$ respectively in $\frak{h}(\mu)$ and $\frak{H}(\mu)$ are
respectively equal to $\sum_{\delta\in\mathcal{D}_{k}^{m}}c_{\delta}^{D_{m}%
}K_{\widehat{\lambda},\widehat{\mu}+\delta}^{A_{m}}$ and $\sum_{\delta
\in\mathcal{D}_{k}^{m}}c_{\delta}^{C_{m}}K_{\widehat{\lambda},\widehat{\mu
}+\delta}^{A_{m}}.$

\item  Consider $\lambda,\mu$ two partitions of length $l$ such that
$k=\frac{\left|  \mu\right|  -\left|  \lambda\right|  }{2}\in\mathbb{N}$ and
set $m=\mathrm{max}(\lambda_{1},\mu_{1})$.\ Then for types $B_{l},C_{l}%
,D_{l},$ the multiplicity of the representation $V(\lambda)$ respectively in
$\frak{e}(\mu)$ and $\frak{E}(\mu)$ are respectively equal to $\sum_{\delta
\in\mathcal{D}_{k}^{m}}c_{\delta}^{D_{m}}K_{\widetilde{\lambda},\widetilde
{\mu}+\delta}^{A_{m}}$ and $\sum_{\delta\in\mathcal{D}_{k}^{m}}c_{\delta
}^{C_{m}}K_{\widetilde{\lambda},\widetilde{\mu}+\delta}^{A_{m}}.$
\end{itemize}
\end{corollary}

\bigskip

\noindent\textbf{Acknowledgments: }The author thanks the organizers of the
workshop ``Combinatorial aspects of integrable systems'' for their hospitality
during the summer 2004 when this work has been completed.\ He would like also
express his gratitude to Professors Okado and Shimozono for many fruitful discussions.

\end{document}